\RequirePackage{fix-cm}

\documentclass[smallextended]{svjour3}       

\smartqed  

\usepackage{graphicx}

\usepackage{latexsym}\usepackage{amsmath}\usepackage{xspace}\usepackage{enumerate}\usepackage{amsfonts}\usepackage{amscd}\usepackage{syntonly}\usepackage{amssymb}

\newtheorem{thm}{Theorem}
\newtheorem{lem}[thm]{Lemma}
\newtheorem{prop}[thm]{Proposition}

\newtheorem{cor}[thm]{Corollary}
\newtheorem{Claim}{Claim}

\newtheorem{rem}[thm]{Remark}

\newenvironment{pf}[1][Proof]{\noindent \textit{#1} }{\hspace{\stretch{1}}}\newenvironment{enui}{\begin{enumerate}[(i)]}{\end{enumerate}}

\newcommand{\BAR}[1]{{\overline{#1}}}

\newcommand\dd{\partial}

\newcommand\ga{\gamma}

\newcommand\eps{\varepsilon}
\newcommand\Om{\Omega}
\newcommand\om{\omega}

\newcommand\Si{\Sigma}
\newcommand\si{\sigma}

\newcommand\HZ{{\operatorname{HZ}}}

\renewcommand\phi{\varphi}

\newcommand\kk{\operatorname{k_\geq}}

\newcommand\keq{\operatorname{k}}
\newcommand{\N}{\mathbb{N}}
\newcommand{\Z}{\mathbb{Z}}

\newcommand{\R}{\mathbb{R}}
\newcommand{\T}{\mathbb{T}}
\newcommand\Cross{{\times}}
\newcommand{\wt}[1]{\widetilde{#1}}
\newcommand\iso{\cong}

\newcommand\inj{\hookrightarrow}
\newcommand\Hd{\dim}
\newcommand\X{\mathcal{X}}

\newcommand\hhat[1]{\widehat{#1}}
\newcommand\nn{\nonumber}

\newcommand\pt{\operatorname{pt}}
\newcommand\K{\operatorname{K}}
\newcommand\C{\mathbb C}
\newcommand\CC{\mathcal{C}}

\newcommand\corank{{\operatorname{corank\,}}}
\newcommand\x{\times}
\newcommand\wo{\setminus}
\newcommand\sub{\subseteq}
\newcommand\one{\mathbf{1}}
\newcommand\Ham{\operatorname{Ham}}

\newcommand\D{\mathbb{D}}

\newcommand\can{{\operatorname{can}}}
\newcommand\id{\operatorname{id}}
\newcommand\pr{\operatorname{pr}}
\newcommand\emb{\operatorname{emb}}
\newcommand\coiso{{\operatorname{coiso}}}
\newcommand\const{\equiv}

\newcommand\cont{{\supseteq}}
\newcommand\HH{\mathcal{H}}

\newcommand\M{\mathcal{M}}
\newcommand\Symp{\operatorname{Symp}}

\journalname{Mathematische Zeitschrift}

\begin{document}

\title{Coisotropic Displacement and Small Subsets of a Symplectic Manifold}

\author{Jan Swoboda\and Fabian Ziltener
}

\institute{Jan Swoboda\at
Max-Planck-Institut f\"ur Mathematik\\
Vivatsgasse 7, D-53111 Bonn, Germany\\
phone: +49 228 402 337\\
\email{swoboda@mpim-bonn.mpg.de}
\and
Fabian Ziltener\at    
Korea Institute for Advanced Study\\
Hoegiro 87 (207-43 Cheongnyangni-dong), Dongdaemun-gu, Seoul 130-722, Korea\\
phone: +82 2 958 3736\\
fabian@kias.re.kr
}

\date{Received: date / Accepted: date}

\maketitle

\begin{abstract}
We prove a coisotropic intersection result and deduce the following:
\begin{itemize}
\item Lower bounds on the displacement energy of a subset of a symplectic manifold, in particular a sharp stable energy-Gromov-width inequality.  
\item A stable non-squeezing result for neighborhoods of products of unit spheres.
\item Existence of a ``badly squeezable'' set in $\R^{2n}$ of Hausdorff dimension at most $d$, for every $n\geq2$ and $d\geq n$.
\item Existence of a stably exotic symplectic form on $\R^{2n}$, for every $n\geq2$.
\item Non-triviality of a new capacity, which is based on the minimal action of a regular coisotropic submanifold of dimension $d$.
\end{itemize}


\end{abstract}

\tableofcontents

\section{Motivation and main results}\label{sec:mot main}
\subsection{Questions}\label{subsec:q} 
The theme of this article is the following.\\

\begin{question}\label{q:small} 
How much symplectic geometry can a small subset of a symplectic manifold carry?
\end{question}
We approach this question from several points of view, interpreting "small" as ``of Hausdorff dimension bounded above by a given number''. Our main tool is an intersection result for coisotropic submanifolds (Theorem \ref{thm:N phi N} below). Further results are proved in \cite{SZHofer}.

One instance of Question \ref{q:small} is the following. Let $(M,\om)$ be a symplectic manifold. (For simplicity all manifolds in this paper are assumed to have empty boundary.) For a subset $X\sub M$ we denote by $e(X):=e(X,M):=e(X,M,\om)$ its displacement energy (see (\ref{eq:e M om X}) below). 
\begin{question}\label{q:e} What lower bounds on $e(X)$ can be detected by special small subsets of $X$? 
\end{question}
To make this more precise, let $\X$ be a collection of subsets of $M$ and $f:\X\to[-\infty,\infty]$ a function. We define 
\[\hhat f:\big\{\textrm{subset of }M\big\}\to[-\infty,\infty],\quad \hhat f(X):=\sup\big\{f(Y)\,\big|\,Y\in\X,\,Y\sub X\big\}.\]
Note that the estimate
\[e(X)\geq \hhat f(X),\quad\forall X\sub M\]
holds, provided that the inequality 
\begin{equation}\label{eq:e X M f}e(X)\geq f(X),\quad\forall X\in\X
\end{equation}
is satisfied. Our goal is therefore to find a collection $\X$ containing a lot of small subsets of $M$ and a big function $f:\X\to[-\infty,\infty]$ for which the inequality (\ref{eq:e X M f}) is satisfied. Our ansatz in this article is to define $\X$ to be the set of all closed regular coisotropic submanifolds $N\sub M$ and $f(N)$ to be a refined version of the minimal symplectic action of $N$ (see (\ref{eq:A Cross}) below). Inequality (\ref{eq:e X M f}) is then a direct consequence of Theorem \ref{thm:N phi N} below.

Another instance of Question \ref{q:small} is the following. Let $n\in\N=\{1,2\ldots\}$ and $a\in(0,\infty)$. We denote by $B^{2n}(a)$ ($\BAR{B}^{2n}(a)$) the open (closed) ball in $\R^{2n}$ of radius $\sqrt{a/\pi}$, around 0. Furthermore, we denote by $Z^{2n}(a):=B^2(a)\x\R^{2n-2}$ the open symplectic cylinder of area $a$. (Note that $B^{2n}(a)$ and $Z^{2n}(a)$ both have Gromov-width $a$.) We abbreviate $B^{2n}:=B^{2n}(\pi)$ and $Z^{2n}:=Z^{2n}(\pi)$. 

Let $(M,\om)$ and $(M',\om')$ be symplectic manifolds. We write $(M,\om)\inj(M',\om')$ iff there exists a symplectic embedding of $M$ into $M'$. Let $n\in\N$. We denote by $\om_0$ the standard symplectic form on $\R^{2n}$. Gromov's non-squeezing result \cite[Corollary, p.~310]{Gr} states that $(B^{2n}(a),\om_0)\not\inj(Z^{2n},\om_0)$ if $a>\pi$. We may ask whether the boundary of the unit ball (or more generally, a finite product of unit spheres) is already too big to be squeezed into the unit cylinder. To make this more precise, let $k\in\N$ and $n_1,\ldots,n_k\in\N$. We denote $n:=\sum_{i=1}^kn_i$. 
\begin{question}[Skinny non-squeezing]\label{q:skinny} If $U\sub\R^{2n}$ is any open neighborhood of $\Cross_{i=1}^kS^{2n_i-1}$ then is it true that
\begin{equation}
(U,\om_0)\not\inj(Z^{2n},\om_0)?
\end{equation}
\end{question}
If $n_i=1$ for some $i$ then an elementary argument shows that there exists a $U$ as above for which $(U,\om_0)\inj(Z^{2n},\om_0)$. Hence assume that $n_i\geq2$, for every $i$. In this case Corollary \ref{cor:skinny} below provides a positive answer to a stabilized version of Question \ref{q:skinny} (and in particular to the original question).

Generalizing Question \ref{q:skinny}, we may wonder how much small subsets of a symplectic manifold can be squeezed. This leads to the following definition. We denote the Hausdorff dimension of a metric space $(X,d)$ by $\Hd(X)$. Let $(M_0,\om_0)$ be a symplectic manifold of dimension $2n$ and $d\in[0,\infty)$. We define the map
\begin{equation}\label{eq:emb M 0 om 0}\emb^{M_0,\om_0}_d:\big\{\textrm{symplectic manifold }(M,\om)\,\big|\,\dim M=2n\big\}\to[0,\infty]\end{equation}
as follows. We define 
\[\X^{2n}:=\big\{(M,\om,X)\,\big|\,(M,\om):\textrm{ symplectic manifold, }\dim M=2n,\,X\sub M\big\},\]
and the map $\emb^{M_0,\om_0}:\X^{2n}\to[0,\infty]$ by 
\begin{eqnarray}
\label{eq:emb M 0 om 0 M om X}&\emb^{M_0,\om_0}(M,\om,X):=&\\
\nn&\inf\big\{a>0\,\big|\,\exists U\sub M\textrm{ open: }X\sub U,\,(U,\om|_U)\inj\big(M_0,a\om_0\big)\big\}.&
\end{eqnarray}
For $d\in[0,\infty)$ we now define the map (\ref{eq:emb M 0 om 0}) by
\[\emb^{M_0,\om_0}_d(M,\om):=\]
\[\sup\big\{\emb^{M_0,\om_0}(M,\om,X)\,\big|\,X\sub M:\textrm{compact},\,\Hd(X)\leq d\big\}.\]
\begin{question}[Squeezing small sets]\label{q:squeezing} What is the value of $\emb^{M_0,\om_0}_d(M,\om)$?
\end{question}
Consider the case $M_0:=Z^{2n}$, equipped with the standard form $\om_0$, and $(M,\om):=(B^{2n},\om_0)$. Then by definition, we have $\emb^{Z^{2n},\om_0}_d(B^{2n},\om_0)\leq\pi$. On the other hand, Theorem \ref{thm:emb d B} below provides a lower bound on this number. In particular, it shows that for every $n\geq2$ and $d\geq n$ the number is positive, and therefore $\emb^{Z^{2n},\om_0}_d$ is an intrinsic symplectic capacity on $(\R^{2n},\om_0)$ in the sense of \cite[Definition C.1, p.~224]{SchlEmb}. 

Looking at Question \ref{q:small} from yet another point of view, we may ask the following.
\begin{question}[Distinguishing symplectic structures]\label{q:dist}How can the \\
(non-)existence of certain subsets be used to distinguish symplectic structures?
\end{question}
Corollary \ref{cor:exotic} below is concerned with this question. It says that every coisotropically infinite symplectic structure on $\R^{2n}$ is stably exotic. (For definitions see page \pageref{stably exotic}.) It follows that there exists a stably exotic symplectic form on $\R^{2n}$, if $n\geq2$. 
\subsection{Coisotropic intersections and displacement energy}\label{subsec:displ coiso}
\subsubsection*{Coisotropic intersections}
The main results of this article are consequences of the following key result. In order to state it, let $(M,\om)$ be a symplectic manifold. We call it \emph{(geometrically) bounded} \label{bounded} iff there exist an almost complex structure $J$ on $M$ and a complete Riemannian metric $g$ such that the following conditions hold:
\begin{itemize}
\item The sectional curvature of $g$ is bounded and $\inf_{x\in M}\iota^g_x>0$, where $\iota^g_x$ denotes the injectivity radius of $g$ at the point $x\in M$. 
\item There exists a constant $C\in(0,\infty)$ such that 
\[|\om(v,w)|\leq C|v|\,|w|,\quad\om(v,Jv)\geq C^{-1}|v|^2,\]
for all $v,w\in T_xM$ and $x\in M$. Here $|v|:=\sqrt{g(v,v)}$.
\end{itemize}
For examples see Section \ref{subsubsec:examples}, page \pageref{ex:bounded}.

Let $N\sub M$ be a coisotropic submanifold. We denote by $A(N)$ its minimal symplectic action (see (\ref{eq:A M om N}) below). We define the \emph{split minimal symplectic action of $N$, $A_\Cross(M,\om,N)$} as follows. We define a \emph{bounded splitting of $(M,\om,N)$} to be a tuple $(M_i,\om_i,N_i)_{i=1,\ldots,k}$, where $k\in\N$ and for every $i=1,\ldots,k$, $(M_i,\om_i)$ is a bounded symplectic manifold and $N_i\sub M_i$ a coisotropic submanifold, such that there exists a symplectomorphism $\phi$ from $\big(\Cross_{i=1}^kM_i,\oplus_{i=1}^k\om_i\big)$ to $(M,\om)$, satisfying $\phi\big(\Cross_{i=1}^kN_i\big)=N$. We define 
\begin{equation}\label{eq:A Cross}A_\Cross(N)=A_\Cross(M,\om,N):=
\end{equation}
\[\sup\big\{\min_{i=1,\ldots k}A(M_i,\om_i,N_i)\,\big|\,(M_i,\om_i,N_i)_i\textrm{ bounded splitting of }(M,\om,N)\big\}.\]
Here our convention is that $\sup\emptyset=0$. 

\noindent{\bf Remark.}~If $(M,\om)$ is not bounded then $(M,\om,N)$ does not admit any bounded splitting, and therefore $A_{\Cross}(N)=0$. This follows from the facts that a finite product of bounded symplectic manifolds is bounded, and boundedness is invariant under symplectomorphisms. $\Box$

We call a coisotropic submanifold $N\sub M$ \emph{regular} iff its isotropy relation (see (\ref{eq:R N om}) below) is a closed subset and a submanifold of $N\x N$. Equivalently, the symplectic quotient of $N$ is well-defined. 

We denote by $\Ham(M,\om)$ the group of Hamiltonian diffeomorphisms on $M$ and by $\Vert\cdot\Vert_\om$ the Hofer norm on $\Ham(M,\om)$. (See Section \ref{sec:background}, page \pageref{Ham M om}). The key result of this article is the following. 

\begin{thm}[Coisotropic intersections]\label{thm:N phi N}Let $(M,\om)$ be a symplectic manifold, $\emptyset\neq N\sub M$ a closed connected regular coisotropic submanifold, and $\phi:M\to M$ a Hamiltonian diffeomorphism. If 
\begin{equation}\label{eq:d phi id}\Vert\phi\Vert_\om<A_\Cross(N)\end{equation} 
then 
\begin{equation}\label{eq:N phi N}N\cap\phi(N)\neq\emptyset.\end{equation}
\end{thm}
In the case where $N$ is a Lagrangian submanifold the statement of this result is an immediate consequence of the Main Theorem in Y.~Chekanov's paper \cite{Ch}. Furthermore, Theorem \ref{thm:N phi N} is related to the main result, Theorem 1, in \cite{ZiLeafwise}. Morally, it has weaker hypotheses, but also a weaker conclusion than that result. (In \cite[Theorem 1]{ZiLeafwise} the condition (\ref{eq:d phi id}) is replaced by a condition involving $A(N)$, which is bounded above by $A_\Cross(N)$. On the other hand, that result provides a lower bound on the number of leafwise fixed points, which is stronger than (\ref{eq:N phi N}).)

The proof of Theorem \ref{thm:N phi N} is an adaption of the proof of \cite[Theorem 1]{ZiLeafwise}. It is based on a certain Lagrangian embedding of $N$ and on the Main Theorem in \cite{Ch}.
\subsubsection*{Displacement energy}
The next result provides an answer to Question \ref{q:e}. To formulate it, we define the map
\[A^\coiso_\Cross:\big\{(M,\om,X)\,|\,(M,\om)\,\textrm{symplectic manifold,}\,X\subseteq M\}\to[0,\infty],\]
by taking $A^{\coiso}_{\Cross}(X)=A^{\coiso}_{\Cross}(M,\om,X)$ to be the supremum of all numbers $A_\Cross(N)$, where $N\neq\emptyset$ is a closed regular coisotropic submanifold of $M$ that is contained in $X$ (with the convention that $\sup\emptyset=0$). Let $(M,\om)$ be a symplectic manifold and $X\sub M$ a subset. 
\begin{cor}[Displacement energy]\label{cor:e A Cross} If $(M,\om)$ is bounded then
\[e(X)\geq A^\coiso_\Cross(X).\]
\end{cor}
This is an immediate consequence of Theorem \ref{thm:N phi N}. In an example on page \pageref{example A Cross} below we will compute $A^\coiso_\Cross(X)$ for certain products. As a special case, let $\om_1$ be an area form on $S^2$ with total area at least $\pi$. Then Corollary \ref{cor:e A Cross}, and inequality (\ref{eq:A min}) and Remarks \ref{rem:e M M'},\ref{rem:e Z 2n a} below imply that
\[e\big((S^3)^2,\R^8,\om_0\big)=\pi,\quad e\big((S^1)^2\x S^2,\R^4\x S^2,\om_0\oplus\om_1\big)=\pi.\]
To our knowledge, these equalities are new. 
\subsubsection*{Stable sharp energy-Gromov-width inequality}
As a consequence of Corollary \ref{cor:e A Cross}, we obtain the following result. Let $(M,\om)$ be a symplectic manifold. We call it \emph{aspherical} iff $\int_{S^2}u^*\om=0$ for every $u\in C^\infty(S^2,M)$. We denote $2n:=\dim M$ and by
\[w(M):=w(M,\om):=\sup\big\{a\,\big|\,(B^{2n}(a),\om_0)\inj(M,\om)\big\}\] 
the Gromov-width of $(M,\om)$. Let $(M',\om')$ be another symplectic manifold.
\begin{cor}[Energy-Gromov-width inequality]\label{cor:e w} Assume that $(M,\om)$ and $(M',\om')$ are aspherical, $(M,\om)$ is bounded, and $M'$ is closed. Then for every open subset $U\sub M$, we have
\begin{equation}\label{eq:e w}e\big(U\x M',M\x M'\big)\geq w(U).
\end{equation}
\end{cor}
This inequality appears to be new. (There are previous results about the case $M'=\{\pt\}$ or with a constant factor on the right hand side of (\ref{eq:e w}), see Section \ref{subsec:related}.)

\noindent{\bf Remark.} Even in the case $M'=\{\pt\}$ the result is sometimes new. As an example, let $X$ be a closed manifold and $\si$ a closed two-form on $X$ such that $\int_{S^2}u^*\si=0$, for every $u\in C^\infty(S^2,X)$. We denote by $\pi:T^*X\to X$ and $\om_\can$ the canonical projection and two-form on $T^*X$. We define $(M,\om):=\big(T^*X,\om_\can+\pi^*\si\big)$ and $M':=\{\pt\}$. 

Then the hypotheses of Corollary \ref{cor:e w} are satisfied, and therefore, applying the corollary, we have $e(U)\geq w(U)$, for every open subset $U\sub M$. For $X$ equal e.g.~to the sphere $S^2$ or the two-torus $\T^2$ this result appears to be new. $\Box$

The statement of the corollary is sharp in the sense that for every pair of symplectic manifolds $(M,\om)$ and $(M',\om')$ there exists an open subset $U\sub M$ for which equality in (\ref{eq:e w}) holds. Namely, denoting $2n:=\dim M$, and by $B^{2n}_r\sub\R^{2n}$ the open ball of radius $r$ around 0, we may choose $r>0$ and an embedding $\phi:(B^{2n}_{3r},\om_0)\inj(M,\om)$. We define $U:=\phi(B^{2n}_r)$. The opposite inequality in (\ref{eq:e w}) then follows from an elementary argument using Remarks \ref{rem:e M M'} and \ref{rem:e Z 2n a} below.
\subsection{The regular coisotropic capacity}
Let $d\in\N$. The minimal coisotropic area gives rise to a map
\[A_\coiso^d:\big\{\textrm{symplectic manifold}\big\}\to[0,\infty].\]
Namely, we define $A_\coiso^d(M,\om)$ to be the supremum of all numbers $A(N)$, where $N\subseteq M$ is a non-empty closed regular coisotropic submanifold of dimension $d$, satisfying the following condition:
\begin{equation}\label{eq:contr}\forall F\textrm{ isotropic leaf of }N,\,\forall x\in C(S^1,F):\,x\textrm{ is contractible in }M.
\end{equation}
Our next result involves the function 
\begin{equation}\label{eq:k}\keq:\N\x\N\to\N\cup\{\infty\},\end{equation}
which is defined as follows. Let $(n,d)\in\N\x\N$. We define $\keq(n,d)$ to be the infimum of all integers $\sum_{i=1}^\ell k_i$, where $\ell\in\N$ and $k_1,\ldots,k_\ell\in\N$ are such that there exist $n_i\in\N$, for $i=1,\ldots,\ell$, satisfying
\begin{eqnarray}\label{eq:n i k i}&n_i\geq k_i,&\\
\label{eq:k i n i d}&\sum_ik_in_i=n,\quad \sum_ik_i(2n_i-k_i)=d.&
\end{eqnarray}
Here our convention is that the infimum of the empty set is $\infty$. Note that $\keq(n,d)=\infty$, if $d<n$ or $d>2n-1$. On the other hand, 
\[\keq(n,d)\leq2n-d,\quad\textrm{if }n\leq d\leq2n-1.\]
(See inequality (\ref{eq:2n d}) in Proposition \ref{prop:k n d} below). Let $n,n'\in\N$, $d\in\{n,\ldots,2n\}$, and $(M,\om)$ and $(M',\om')$ be symplectic manifolds of dimensions $2n$ and $2n'$, respectively. 
\begin{thm}[Regular coisotropic capacity]\label{thm:A coiso cap} The following statements hold. 
\begin{enui}
\item\label{thm:A coiso cap:cap} If $d<2n$ then the restriction of $A_{\coiso}^d$ to the class of all aspherical symplectic manifolds of dimension $2n$ is a symplectic capacity (as defined on page \pageref{cap}). For $d=2n-1$ this capacity is normalized, i.e., it takes on the value $\pi$ on $B^{2n}$ and $Z^{2n}$. 
\item\label{thm:A coiso cap:A coiso d B Z} We have
\begin{eqnarray}\label{eq:A coiso d B}&A_\coiso^d(B^{2n},\om_0)\geq\frac\pi{\keq(n,d)},&\\
\label{eq:A coiso d Z}&A_\coiso^d(Z^{2n},\om_0)\leq\pi.&
\end{eqnarray}
\item\label{thm:A coiso cap:asph} If $(M',\om')$ is closed and aspherical then 
\begin{equation}\label{eq:A coiso d 2n'}A_\coiso^{d+2n'}\big(M\x M',\om\oplus\om'\big)\geq A_\coiso^d(M,\om).
\end{equation}
\end{enui}
\end{thm}
For $d\in\{n,\ldots,2n-1\}$ we call the restriction of $A_\coiso^d$ to the class of all aspherical symplectic manifolds of dimension $2n$ the {\emph{regular coisotropic capacity}}\label{coiso cap}. In the case $d=n$ this is closely related to the Lagrangian capacity introduced by K.~Cieliebak and K.~Mohnke \cite{CM}. Furthermore, in this case, the right hand side in (\ref{eq:A coiso d B}) simplifies. Namely, we have 
\begin{equation}\keq(n,n)=\K(n):=\inf\big\{\sum_{i=1}^\ell k_i\,\big|\,\ell\in\N,\,k_1,\ldots,k_\ell\in\N:\,n=\sum_ik_i^2\big\}.\label{eq:k n inf}
\end{equation}
(See equality (\ref{eq:k n d:k n n}) in Proposition \ref{prop:k n d} below.) The first few values of $\K$ are 
\[\begin{array}{rrrrrrrrrrrrrrrrrrrr}
    n=1&2&3&4&5&6&7&8&9&10&11&12&13&14&15&16&17&18&19&20\\
\K(n)=1&2&3&2&3&4&5&4&3& 4& 5& 6& 5& 6& 7& 4& 5& 6 & 7& 6
\end{array}
\]
The function $\K$ satisfies the upper bound 
\[\K(n)<\sqrt n+2^{\frac32}\sqrt[4]n.\]
(See inequality (\ref{eq:k n d:sqrt n}) in Proposition \ref{prop:k n d} below.) Via inequality (\ref{eq:A coiso d B}) this yields a lower bound on $A_\coiso^d(B^{2n},\om_0)$. 
\subsection{Symplectic squeezing of small sets}\label{subsec:squeezing}
Our next application of Theorem \ref{thm:N phi N} is the following. Let $k\in\N$ and $n_1,\ldots,n_k\in\N$. We denote $n:=\sum_{i=1}^kn_i$. Let $(M,\om)$ be a symplectic manifold. 
\begin{cor}[Skinny non-squeezing]\label{cor:skinny} Let $U\sub\R^{2n}$ be an open neighborhood of $\Cross_{i=1}^kS^{2n_i-1}$. Assume that $n_i\geq2$, for every $i=1,\ldots,k$, $M$ is closed and connected, and 
\begin{equation}\label{eq:int u om}\int_{S^2}u^*\om\in\pi\Z,\quad\forall u\in C^\infty(S^2,M).\end{equation} 
Then we have
\begin{equation}\label{eq:U om 0}\big(U\x M,\om_0\oplus\om\big)\not\inj\big(Z^{2n}\x M,\om_0\oplus\om\big).
\end{equation}
\end{cor}
Taking $M=\{\pt\}$, this result provides a positive answer to Question \ref{q:skinny} of Section \ref{subsec:q} in the case $n_i\geq2$, for every $i$. To our knowledge, this fact is new.\footnote{We are not aware of any other written proof. \label{skinny neck}However, K.~Cieliebak mentioned to the second author that in the case $M=\{\pt\}$ the non-embedding (\ref{eq:U om 0}) also follows from a standard ``neck stretching'' argument similar to the proof of \cite[Theorem 1.1]{CM}.}  

In order to state our result about Question \ref{q:squeezing}, recall the definition (\ref{eq:emb M 0 om 0}) of the map $\emb^{M_0,\om_0}_d$. We also need the following. We define the map 
\begin{equation}\label{eq:k N 0}\kk:\N\x[0,\infty)\to\R\end{equation}
in the same way as the map $\keq$ (see (\ref{eq:k})), replacing (\ref{eq:k i n i d}) by the conditions 
\begin{equation}
\label{eq:k i n i d geq}\sum_ik_in_i\geq n,\quad \sum_ik_i(2n_i-k_i)\leq d.
\end{equation}
Note that 
\[\kk(n,d)\leq\keq(n,d),\quad\forall(n,d)\in\N\x[0,\infty).\]
We can now formulate the following result. 
\begin{thm}[Badly squeezable small sets]\label{thm:emb d B} For every $n\in\{2,3,\ldots\}$ and $d\in[n,\infty)$ we have 
\[\emb^{Z^{2n},\om_0}_d(B^{2n},\om_0)\geq \frac\pi{\kk(n,d)}.\]
\end{thm}
The map $\kk$ satisfies some explicit upper bounds, see Proposition \ref{prop:k n d} below. As an example, inequalities (\ref{eq:k d d'},\ref{eq:k n d:k n n},\ref{eq:k n d:sqrt n}) of that proposition imply that for $d\geq n$, 
\[\kk(n,d)<\sqrt n+2^{\frac32}\sqrt[4]n.\]
In the proof of Theorem \ref{thm:emb d B} we consider a certain product of Stiefel manifolds. This is a regular coisotropic submanifold $N$ of $\R^{2(n+n')}$ for some $n'\in\N$. We also use the inequality 
\[\emb^{Z^{2n,\om_0}}_d(M,\om)\geq A_\coiso^d(M,\om),\]
see Proposition \ref{prop:ineqembA} below. The proof of this inequality is based on Theorem \ref{thm:N phi N}. It also relies on an argument in which we glue disks to a given regular coisotropic submanifold $N$, to make all loops in the isotropic fibers of $N$ contractible. 
\subsection{Stably exotic symplectic forms}\label{subsec:exotic}
To state our last application of Theorem \ref{thm:N phi N}, let $n\in\N$. We call a symplectic form $\om$ on $\R^{2n}$ \emph{stably exotic} \label{stably exotic} iff the following holds. Let $(X,\si,\si')$ be a triple consisting of a closed manifold $X$ and symplectic forms $\si$ and $\si'$ on $X$, with $\si$ aspherical, and let $\phi:M:=X\x\R^{2n}\to X\x\R^{2n}$ be an embedding. Then 
\begin{equation}\label{eq:phi si om 0}\phi^*(\si\oplus\om_0)\neq\si'\oplus\om.
\end{equation}
Note that such an $\om$ is exotic in the usual sense, i.e., if $\phi:\R^{2n}\to\R^{2n}$ is an embedding then $\phi^*\om_0\neq\om$. 

Our result is a sufficient criterion for stable exoticness. Namely, we call a symplectic manifold $(M,\om)$ \emph{(coisotropically) infinite} iff there exists a non-empty regular closed coisotropic submanifold $N\sub M$ with $A(N)=\infty$. (This means that $\int_{S^2}u^*\om=0$, for every $u\in C^\infty(S^2,M)$ such that $u(S^1)$ is contained in some isotropic leaf of $N$.) 
\begin{cor}[Stably exotic form]\label{cor:exotic} Let $n\in\N$. Then every coisotropically infinite form on $\R^{2n}$ is stably exotic.
\end{cor}
It follows from Corollary \ref{cor:exotic} and 
an example on page \pageref{examples infinite} below that for every $n\geq2$, there exists a stably exotic symplectic form on $\R^{2n}$.
\subsection{Remarks, examples, related work, organization of the article}\label{subsec:related}
\subsubsection*{Remarks}\label{subsubsec:remarks}
\noindent{\bf On geometric boundedness:}~In the article \cite{ZiLeafwise} the second author used the term ``geometrically bounded'' in a slightly stronger sense. 

\noindent{\bf On Theorem \ref{thm:N phi N}:}~Assume that the hypotheses of Theorem \ref{thm:N phi N} are satisfied and the pair $(N,\phi)$ is non-degenerate in the sense of \cite{ZiLeafwise}. Then the number of leafwise fixed points of $\phi$ is bounded below by the sum of the $\Z_2$-Betti numbers of $N$. This follows by adapting the proof of Theorem \ref{thm:N phi N} along the lines of the proof of \cite[Theorem 1]{ZiLeafwise}. $\Box$

\noindent{\bf On Corollary \ref{cor:e w}:} The assumption that $(M',\om')$ is aspherical can be weakened as follows. For a symplectic manifold $(M,\om)$ we define
\begin{equation}\label{eq:A M om}A(M,\om):=\inf\big(\big\{\int_{S^2}u^*\om\,\big|\,u\in C^\infty(S^2,M)\big\}\cap(0,\infty)\big)\in[0,\infty].\end{equation}
Assume that $(M,\om)$ is aspherical and bounded, and there exist closed symplectic manifolds $(M_i,\om_i)$, for $i=1,\ldots,k$, such that $M'=\prod_iM_i$, $\om'=\oplus_i\om_i$, and $A(M_i,\om_i)\geq w(U)$. Then inequality (\ref{eq:e w}) still holds. This follows from an argument using Corollary \ref{cor:e A Cross}. $\Box$

\noindent{\bf On the regular coisotropic capacity:} One can define variants of this capacity by imposing other conditions on the coisotropic submanifold $N$ (e.g., stability or a contact type condition). Note that in order to obtain a capacity $c^d$ satisfying $c^d(Z^{2n},\om_0)<\infty$, one can neither completely drop the condition that $N$ is closed nor that it is regular.

Namely, there exists a regular (but not closed) coisotropic submanifold $N\sub Z^{2n}$ such that $A(Z^{2n},\om_0,N)=\infty$, and there exists a closed (but not regular) coisotropic submanifold $N'\sub Z^{2n}$ such that $A(Z^{2n},\om_0,N')=\infty$. 

As an example, we may choose a coisotropic subspace $W\sub\R^{2n}$ of dimension $d$ and define $N:=W\cap Z^{2n}$. Furthermore, we may choose $N'$ to be a closed hypersurface in $Z^{2n}$ without any closed characteristic. Such an $N'$ exists by a construction due to V.~Ginzburg, see \cite{Gi}, Example 7.2 p.~158. (We shrink Ginzburg's hypersurface homothetically, so that it fits into $Z^{2n}$. Such an $N$ is not regular.) $\Box$

\subsubsection*{Examples}\label{subsubsec:examples}
\noindent{\bf Examples for geometric boundedness:}\label{ex:bounded} $(M,\om)$ is bounded if it is closed (i.e., compact and with empty boundary), a symplectic vector space, convex at infinity (see \cite[Remark 2.3]{CGK}), or the twisted cotangent bundle $\big(T^*X,\om_\can+\pi^*\si\big)$, where $X$ is a closed manifold $X$, $\si$ is a closed two-form on $X$, and $\pi:T^*X\to X$ and $\om_\can$ denote the canonical projection and two-form on $T^*X$. (For the last example see \cite[Proposition 2.2]{CGK}.) Furthermore, by a straight-forward argument, the product of a two bounded symplectic manifolds is bounded. $\Box$

\noindent{\bf Example for $A^\coiso_\Cross(X)$ and Corollary \ref{cor:e A Cross}:}\label{example A Cross}
For $n\in\N$ and $a\in(0,\infty)$ we denote by $S^{2n-1}(a)\sub\R^{2n}$ the sphere of radius $\sqrt{a/\pi}$ around 0. Let $k,\ell\in\N_0=\{0,1,\ldots\}$, for $i=1,\ldots,k$ let $n_i\in\N$ and $a_i\in(0,\infty)$, and for $i=1,\ldots,\ell$ let $(M_i,\om_i)$ be a closed symplectic manifold. We define $n:=\sum_in_i$, and  
\begin{eqnarray*}&M:=\R^{2n}\x\Cross_{i=1}^\ell M_i,\quad \om:=\om_0\oplus\oplus_{i=1}^\ell\om_i,&\\
&X:=\Cross_{i=1}^kS^{2n_i-1}(a_i)\x\Cross_{i=1}^\ell M_i.&
\end{eqnarray*}
(If $k=0$ or $\ell=0$ then our convention is that the corresponding Cartesian product is a singleton.) We claim that 
\begin{equation}\label{eq:A min} A^\coiso_\Cross(X)\geq a:=\inf\big(\{a_i\}_i\cup\{A(M_i,\om_i)\}_i\big),
\end{equation}
where $A(M_i,\om_i)$ is defined as in (\ref{eq:A M om}). (Here our convention is that $\inf\emptyset=\infty$.) To see this, observe that $N:=X$ is a closed regular coisotropic submanifold of $\R^{2n}\x\Cross_{i=1}^\ell M_i$ with $A_\Cross(N)\geq a$. (This inequality follows from a straight-forward argument involving the splitting $M=\Cross_i\R^{2n_i}\x\Cross_iM_i$, and Remark \ref{rmk:S N N'} and Proposition \ref{minimalareaStiefel} below.) The claimed inequality (\ref{eq:A min}) follows. 

Combining Corollary \ref{cor:e A Cross} with inequality (\ref{eq:A min}) we obtain 
\begin{equation}\label{eq:e a}e(X)\geq a.
\end{equation}
Assume that $k\geq1$ and $\min_{i=1,\ldots,k}a_i\leq\inf_{i=1,\ldots,\ell}A(M_i,\om_i)$. Then the estimate (\ref{eq:e a}) is sharp. To see this, we choose $j\in\{1,\ldots,k\}$ such that $a_j=\min_ia_i$. Remarks \ref{rem:e M M'} and \ref{rem:e Z 2n a} below imply that $e(X)\leq a_j=a$. $\Box$

\noindent{\bf Examples of coisotropically infinite manifolds.}\label{examples infinite} Every closed aspherical symplectic manifold and the cotangent bundle of a closed manifold are infinite. (In the first example we may take $N:=M$, and in the second example we may take $N$ to be the zero section of the bundle.) Furthermore, by a standard argument, there exists a pair $(\om,L)$, where $\om$ is a symplectic form on $\R^{2n}$ and $L\sub\R^{2n}$ is a closed Lagrangian submanifold such that $A(L)=\infty$. (See e.g.~\cite{ALP}, p.~317.) Such a form $\om$ is coisotropically infinite. $\Box$
\subsubsection*{Related work}
\noindent{\bf Coisotropic intersections and displacement energy.} In \cite{ZiLeafwise} the second author proved a result (Theorem 1) similar to the key result (Theorem \ref{thm:N phi N}) of the present article. 

Let $(M,\om)$ be a symplectic manifold and $U\sub M$ an open subset. In the case $M'=\{\pt\}$ the energy-Gromov-width inequality (\ref{eq:e w}) follows from an elementary argument, whenever one can prove that a certain symplectic capacity $c$ satisfies 
\begin{eqnarray}\label{eq:c B}&c(B^{2n},\om_0)\geq\pi,&\\
\label{eq:c U om U}&c(U,\om|_U)\leq e(U,M).&
\end{eqnarray}
In the following, we take ``$e$'' in inequality (\ref{eq:c U om U}) to mean variants of the displacement energy. Let $(M,\om):=(\R^{2n},\om_0)$ and $c$ be the {\bf Ekeland-Hofer capacity}. Then inequality (\ref{eq:c B}) was proved by I.~Ekeland and H.~Hofer in \cite[Theorem 1]{EH}. (They actually proved equality.) Furthermore, H.~Hofer \cite[Theorem 1.6(i)]{HoTop} proved inequality (\ref{eq:c U om U}). 

Let now $c$ be the {\bf $\pi_1$-sensitive Hofer-Zehnder capacity} $c^\circ_\HZ$. Then inequality (\ref{eq:c B}) is an easy consequence of the definition. Furthermore, H.~Hofer \cite[Theorem 2]{HoEst} proved inequality (\ref{eq:c U om U}) for $(\R^{2n},\om_0)$. (See also H.~Hofer and E.~Zehnder \cite[Section 5.5]{HZ}.) 

U.~Frauenfelder, V.~Ginzburg, and F.~Schlenk \cite[Corollary 1]{FGS} proved the inequality (\ref{eq:c U om U}) (for $c=c^\circ_\HZ$) if an exhaustion of $(M,\om)$ admits an action selector. As an example, assume that $(M,\om)$ is aspherical. If it is also closed or convex at infinity, then it admits such an exhaustion. (It even admits an action selector itself. See the examples in \cite[pages 3,4]{FGS}. See also inequality (2.9), p.~13, and Proposition 3.4 in \cite{Gi}.) 

M.~Usher \cite[Corollary 1.2]{Us} proved that (\ref{eq:c U om U}) holds if $(M,\om)$ is of type (C) (see \cite[p.~3]{Us}). Examples of type (C) manifolds are Stein manifolds, closed manifolds, and convex symplectic 4-manifolds. Non-sharp versions of inequality (\ref{eq:c U om U}) were proved by M.~Schwarz \cite[Corollary 5.16]{Schw} and F.~Schlenk \cite[Theorem 1.1]{SchlApp}. 

In \cite{LMGeo} F.~Lalonde and D.~McDuff proved that $e(U,M)\geq\frac12w(U)$, for \emph{any} symplectic manifold $(M,\om)$ (and any open subset $U\sub M$). (This inequality is a non-sharp version of (\ref{eq:e w}) with $M'=\{\pt\}$.) F.~Lalonde and C.~Pestieau \cite{LP} proved the following stabilized version of this result:
Let $(M,\om)$ and $(M',\om')$ be symplectic manifolds, with $M'$ closed, and $U\sub M$ an open subset. Then Theorem 1.2 in \cite{LP} states that 
\[e\big(U\x M',M\x M'\big)\geq \frac12w(U).\] 

\noindent{\bf The regular coisotropic capacity.} Let $n\in\N$. We denote
\begin{align*}\M:=\big\{&(M,\om)\textrm{ symplectic manifold }\big|\\
&\dim M=2n,\,\pi_1(M)\iso\pi_2(M)\iso\{e\}\big\}.\end{align*}
In \cite{CM} K.~Cieliebak and K.~Mohnke defined the \emph{Lagrangian capacity} to be the map 
\[c_L:\M\to[0,\infty),\]
\[c_L(M,\om):=\sup\big\{A(M,\om,L)\,\big|\,L\sub M\textrm{ embedded Lagrangian torus}\big\}.\]
(See also \cite{CHLS}, Sec.~2.4, p.~11.) The authors proved that 
\begin{equation}\label{eq:c L}c_L(B^{2n},\om_0)=\frac\pi n.\end{equation} 
The capacity $c_L$ is bounded above by the regular coisotropic capacity $A_\coiso^n$, since every Lagrangian submanifold is regular. Let $d\in\{n,\ldots,2n-1\}$ be an integer. If $d=n$ then assume that $n\geq4$. Then we have
\begin{equation}\label{eq:c L A d}c_L(B^{2n},\om_0)<A_\coiso^d(B^{2n},\om_0).\end{equation}
To see this, observe that $n>\keq(n,d)$. (In the case $d=n$ this follows by taking $\ell:=n-3$, $k_1:=2$, $k_2,\ldots,k_{n-3}:=1$, and in the case $d>n$ from inequality (\ref{eq:2n d}) below.) Combining this with the inequality (\ref{eq:A coiso d B}) of Theorem \ref{thm:A coiso cap}(\ref{thm:A coiso cap:A coiso d B Z}) and the equality (\ref{eq:c L}), inequality (\ref{eq:c L A d}) follows. 

\noindent{\bf Squeezing small sets.} As mentioned on page \pageref{skinny neck}, an argument by K.~Cieliebak and K.~Mohnke as in \cite{CM} yields the statement of Corollary \ref{cor:skinny} in the case $M=\{\pt\}$. 

\noindent{\bf Exotic symplectic structures.} In \cite{Gr} M.~Gromov proved that there does not exist a closed exact $\om_0$-Lagrangian submanifold of $\R^{2n}$. It was folklore that this property of $(\R^{2n},\om_0)$ implies the existence of an exotic symplectic structure on $\R^{2n}$, and a proof of this appeared in the paper \cite{Vi} by C.~Viterbo. Another reference is \cite[p.~317]{ALP}.

\subsubsection*{Organization of the article}

In Section \ref{sec:background} we collect some definitions that are used throughout this article. We also prove some useful properties of the functions $\K,\keq$, and $\kk$, including upper bounds (Proposition \ref{prop:k n d}). Section \ref{sec:proofs} is devoted to the proofs of the results of Section \ref{sec:mot main}. Since Theorem \ref{thm:A coiso cap} is used in the proof of Corollary \ref{cor:e w}, we prove it before the corollary. Appendix \ref{sec:aux} contains some basic facts from (pre-)symplectic geometry, topology, and manifold theory, which are used in the proofs of the main results. 
\section{Background and a further result}\label{sec:background}

In this section some standard symplectic geometry is recalled, which is used in this article. We also prove upper estimates on the functions $\keq$ and $\kk$.

Let $(M,\om)$ be a symplectic manifold. We define the group $\Ham(M,\om)$ of \emph{Hamiltonian diffeomorphisms of $M$}\label{Ham M om}, as follows. We define $\HH(M,\om)$\label{HH M om} to be the set of all functions $H\in C^\infty\big([0,1]\x M,\R\big)$ whose Hamiltonian time-$t$ flow $\phi_H^t:M\to M$ exists and is surjective, for every $t\in[0,1]$. We define 
\[\Ham(M,\om):=\big\{\phi_H^1\,\big|\,H\in\HH(M,\om)\big\}.\]
This is a subgroup of the group of diffeomorphisms of $M$. (See for example \cite{SZHofer}.) It contains the group $\Ham_c(M,\om)$ of Hamiltonian diffeomorphisms generated by a compactly supported time-dependent function. By definition, the Hofer norm on the space of functions is the map 
\[\Vert\cdot\Vert:C^\infty([0,1]\x M,\R)\to[0,\infty],\]
\[\Vert H\Vert:=\int_0^1\big(\sup_MH^t-\inf_MH^t\big)dt,\]
where $H^t(x):=H(t,x)$. (It follows from Lemma \ref{le:Borel} below that this norm is well-defined.) We define the Hofer norm on $\Ham(M,\om)$ to be the map 
\[\Vert\cdot\Vert_\om:\Ham(M,\om)\to[0,\infty],\]
\begin{equation}\label{eq:Vert phi}\Vert\phi\Vert_\om:=\inf\big\{\Vert H\Vert\,\big|\,H\in\HH(M,\om):\,\phi_H^1=\phi\big\}.
\end{equation}
We define the \emph{displacement energy} of a subset $X\sub M$ to be 
\begin{equation}\label{eq:e M om X}e(X,M):=e(X,M,\om):=\inf\big\{\Vert H\Vert\,\big|\,H\in\HH(M,\om)\,\big|\,\phi_H^1(X)\cap X=\emptyset\big\}.
\end{equation}
Let $(M,\om)$ be a symplectic manifold and $N\sub M$ a submanifold. Then $N$ is called \emph{coisotropic} iff for every $x\in N$ the subspace
\[T_xN^\om=\big\{v\in T_xM\,\big|\,\om(v,w)=0,\,\forall w\in T_xN\big\}\]
of $T_xM$ is contained in $T_xN$. As an example, every hypersurface in $M$ is coisotropic.

Let $N\sub M$ be a coisotropic submanifold. We define the \emph{isotropy relation} to be the set 
\begin{align}\label{eq:R N om} R^{N,\om}:=\big\{&(x(0),x(1))\,\big|\,x\in C^\infty([0,1],N):\\
\nn&\dot x(t)\in(T_{x(t)}N)^\om,\,\forall t\in[0,1],\,x(i)=x_i,\,\forall i=0,1\big\}.
\end{align}
This is an equivalence relation on $N$. For a point $x_0\in N$ we call the $R^{N,\om}$-equivalence class of $x_0$ the \emph{isotropic leaf} through $x_0$. We denote this subset of $N$ by $N^\om_{x_0}$. Furthermore, we denote
\[N_\om:=\big\{\textrm{isotropic leaf of }N\big\}.\] 

We call $N$ \emph{regular} \label{regular} if $R^{N,\om}$ is a closed subset and a submanifold of $N\times N$. This holds if and only if there exists a manifold structure on the set $N_\om$ such that the canonical projection $\pi_N:N\to N_{\om}$ is a submersion, cf.~\cite[Lemma 15]{ZiLeafwise}. If $N$ is closed then by C.~Ehresmann's theorem this implies that $\pi_N$ is a smooth (locally trivial) fiber bundle. (See the proposition on p.~31 in \cite{Eh}.) 

We define the \emph{action (or area) spectrum} and the \emph{minimal action} of $N$ as
\begin{equation}\label{eq:S M om N}S(M,\om,N):=\left\{\int_\D u^*\om\,\bigg|\,u\in C^\infty(\D,M):\,\exists F\in N_\om:\, u(S^1)\sub F\right\},
\end{equation}
\begin{equation}\label{eq:A M om N} A(N)=A(M,\om,N):=\inf\big(S(N)\cap (0,\infty)\big)\in[0,\infty].
\end{equation}

Let $n\in\N$. We denote by $\Symp^{2n}$ the class of all symplectic manifolds of dimension $2n$. Let $\CC\sub\Symp^{2n}$ be a subclass with the following properties. We have $(B^{2n},\om_0),(Z^{2n},\om_0)\in\CC$. Furthermore, if $(M,\om)\in\CC$ and $a\in\R\wo\{0\}$, then $(M,a\om)\in\CC$. 

By a \emph{symplectic capacity on $\CC$}\label{cap} we mean a map $c:\CC\to [0,\infty]$, such that for every $(M,\om),(M',\om')\in\CC$, the following conditions are satisfied:\\

{\bf(Monotonicity)} $c(M,\om)\leq c(M',\om')$, if $(M,\om)\inj(M',\om')$,\\

{\bf(Conformality)} $c(M,a\om)=|a|c(M,\om)$, for every $a\in\R\wo\{0\}$,\\

{\bf(Nontriviality)} $0<c(B^{2n})$ and $c(Z^{2n})<\infty$.\\

The next result summarizes some properties of the functions $\K$, $\keq$, $\kk$ (see (\ref{eq:k n inf},\ref{eq:k},\ref{eq:k N 0})). In particular, it provides upper bounds on these functions. 
\begin{prop}\label{prop:k n d} Let $n,n',k\in\N$ and $d,d'\in[0,\infty)$. Then the following \\
(in-)equalities hold:
\begin{eqnarray}\label{eq:k d d'}&\kk(n,d)\geq\kk(n,d'),\quad\textrm{if }d\leq d',&\\
\label{eq:keq kk}&\keq(n,d)\geq\kk(n,d),&\\
\label{eq:k n d:k n n}&\keq(n,n)=\kk(n,n)=\K(n),&\\
\label{eq:k n d:sqrt n}&\K(n)<\sqrt n+2^{\frac32}\sqrt[4]n,&\\ 
\label{eq:2n d}&\keq(n,d)\leq 2n-d,\quad\textrm{if }n\leq d\leq2n-1&\\
\label{eq:sqrt 2n d}&\kk(n,d)<\sqrt{2n-d}+3,\quad\textrm{if }n\geq9,\,n+6\sqrt n-9\leq d\leq2n,&\\
\label{eq:k n n' d d'}&\kk(n+n',d+d')\leq\kk(n,d)+\kk(n',d'),&\\
\label{eq:k n}&\keq(n,2n-k^2)\leq k,\quad\textrm{if }k\textrm{ divides }n\textrm{ and }k^2\leq n.&
\end{eqnarray}
\end{prop}
For the proof of Proposition \ref{prop:k n d}, we need the following.
\begin{rem}\label{rem:max} \upshape
For every $m\in\N$ and $a>0$ we have
\[\max\big\{\sum_{i=1}^mx_i\,\big|\,x_i\in\R,\,\forall i=1,\ldots,m,\,\sum_ix_i^2=a\big\}=\sqrt{ma}.\]
(The maximum is attained at the point $\sqrt{\frac am}(1,\ldots,1)$.) $\Box$
\end{rem}
\begin{proof}[of Proposition \ref{prop:k n d}]\setcounter{Claim}{0} {\bf Inequalities (\ref{eq:k d d'},\ref{eq:keq kk})} are direct consequences of the definitions. 

We show that the {\bf equalities (\ref{eq:k n d:k n n})} hold: We claim that
\begin{equation}\label{kk n n geq K n}\kk(n,n)\geq\K(n).
\end{equation}
To see this, let $\ell\in\N$ and $k_1,\ldots,k_\ell$ be as in the definition of $\kk(n,d)$ with $d=n$. Inequality (\ref{kk n n geq K n}) is a consequence of the next claim.
\begin{Claim}\label{claim:sum k i 2 n} We have
\begin{equation}\label{eq:sum k i 2 n}\sum_{i=1}^\ell k_i^2=n.
\end{equation}
\end{Claim}
\begin{proof}[of Claim \ref{claim:sum k i 2 n}] We choose integers $n_1,\ldots,n_\ell$ such that the inequalities (\ref{eq:n i k i},\ref{eq:k i n i d geq}) are satisfied. Subtracting the first from the second inequality in (\ref{eq:k i n i d geq}), we obtain $\sum_ik_i(n_i-k_i)\leq d-n=0$. Using the inequalities (\ref{eq:n i k i}), it follows that $n_i=k_i$, for every $i=1,\ldots,\ell$. Combining this with (\ref{eq:k i n i d geq}), the equality (\ref{eq:sum k i 2 n}) follows. This proves Claim \ref{claim:sum k i 2 n}.
\end{proof}

We complete the proof of (\ref{eq:k n d:k n n}): In view of (\ref{kk n n geq K n}) and the inequality (\ref{eq:keq kk}) it suffices to show that 
\begin{equation}\label{eq:keq n n K}\keq(n,n)\leq\K(n).\end{equation} 
To see this, let $\ell\in\N$ and $k_1,\ldots,k_\ell\in\N$ be as in the definition of $\K(n)$. This means that $\sum_{i=1}^\ell k_i^2=n$. We define $n_i:=k_i$, for $i=1,\ldots,\ell$. Then the conditions (\ref{eq:n i k i},\ref{eq:k i n i d}) in the definition of $\keq(n,n)$ are satisfied with $d=n$. Inequality (\ref{eq:keq n n K}) follows. This proves (\ref{eq:k n d:k n n}).

To prove {\bf inequality (\ref{eq:k n d:sqrt n})}, let $n\in\N$. We define $\ell:=5$ and $k_1$ to be the biggest integer $\leq\sqrt n$. By the Four Squares Theorem there exist integers $k_2,\ldots,k_5\in\N_0=\{0,1,\ldots\}$ such that $\sum_{i=2}^5k_i^2=n-k_1^2$. (See for example Theorem 2.10 in the book \cite{EW}.) Since $n=\sum_{i=1}^5k_i^2$, by definition, we have
\begin{equation}\label{eq:k n d sum}\K(n)\leq\sum_{i=1}^5k_i.
\end{equation}
Furthermore, by Remark \ref{rem:max} with $m:=4$ and $a:=n-k_1^2$, we have $\sum_{i=2}^5k_i\leq2\sqrt{n-k_1^2}.$ Combining this with the inequality $k_1>\sqrt n-1$, we obtain 
\[\sum_{i=1}^5k_i<k_1+2\sqrt{n-(\sqrt n-1)^2}<\sqrt n+2^{\frac32}\sqrt[4]n.\]
Combining this with (\ref{eq:k n d sum}), inequality (\ref{eq:k n d:sqrt n}) follows. 

{\bf Inequality (\ref{eq:2n d})} follows by taking $\ell:=2n-d$, $k_i:=1$, for $i=1,\ldots,\ell$, $n_i:=1$, for $i=1,\ldots,\ell-1$, and $n_\ell:=d-n+1$.

To show {\bf (\ref{eq:sqrt 2n d})}, assume that $n\geq9$ and $n+6\sqrt n-9\leq d\leq2n$. For every number $x\in\R$ we denote by $\lceil x\rceil$ the smallest integer $\geq x$. We define $\ell:=1$, $k_1:=\lceil\sqrt{2n-d}\rceil+2$, and $n_1:=\lceil\frac n{k_1}\rceil$. The claimed inequality is now a consequence of the following claim. 
\begin{Claim}\label{claim:n i k i d} The conditions (\ref{eq:n i k i},\ref{eq:k i n i d geq}) are satisfied.
\end{Claim}
\begin{proof}[of Claim \ref{claim:n i k i d}] We prove that {\bf condition (\ref{eq:n i k i})} holds. The assumption $d\geq n+6\sqrt n-9$ implies that $2n-d\leq n-6\sqrt n+9=(\sqrt n-3)^2$. Since $n\geq9$, we have $\sqrt n-3\geq0$. It follows that 
\begin{equation}\label{eq:n sqrt}(\sqrt{2n-d}+3)^2\leq n.\end{equation}
On the other hand, we have $k_1<\sqrt{2n-d}+3$, and therefore $n_1\geq\frac n{k_1}>\frac n{\sqrt{2n-d}+3}$. Combining this with (\ref{eq:n sqrt}), it follows that $n_1>k_1$. This proves condition (\ref{eq:n i k i}).

The {\bf first condition in (\ref{eq:k i n i d geq})}, $k_1n_1\geq n$, follows from the definition of $n_1$. 

To prove the {\bf second condition in (\ref{eq:k i n i d geq})}, observe that $2n_1-k_1<\frac{2n}{k_1}+2-\sqrt{2n-d}-2$, and therefore 
\[k_1(2n_1-k_1)<2n-k_1\sqrt{2n-d}\leq2n-(2n-d)=d.\]
This proves the second condition in (\ref{eq:k i n i d geq}), and completes the proof of Claim \ref{claim:n i k i d}, and hence of (\ref{eq:sqrt 2n d}). 
\end{proof}
{\bf Inequality (\ref{eq:k n n' d d'})} follows from a straight-forward argument, and {\bf inequality (\ref{eq:k n})} follows by choosing $\ell:=1$, $k_1:=k$, and $n_1:=n/k$. 

This completes the proof of Proposition \ref{prop:k n d}.
\end{proof}
\section{Proofs of the main results}\label{sec:proofs}
\subsection{Proof of Theorem \ref{thm:N phi N} (Coisotropic intersections)}\label{subsec:proof:thm:N phi N}
A central ingredient of the proof of Theorem \ref{thm:N phi N} is the following result by Y.~Chekanov. Let $(M,\om)$ be a symplectic manifold, $\Si$ a Riemann surface, and $X\sub M$ a subset. For every almost complex structure $J$ on $M$, we define
\begin{eqnarray}\label{eq:A Si M J X}&A\big(\Si,M,\om,J,X\big):=&\\
\nn&\inf\left(\left\{\int_\Si u^*\om\,\big|\,u:\Si\to M\,J\textrm{-holomorphic, }u(\dd\Si)\sub X\right\}\cap(0,\infty)\right).&
\end{eqnarray}
Furthermore, we define the \emph{bounded minimal action of $(M,\om)$ relative to $X$} to be
\begin{equation}\label{eq:A b}
A_b(X)=A_b(M,\om,X):=\end{equation}
\[\sup\big\{\min\{A(S^2,M,\om,J,\emptyset),A(\D,M,\om,J,X)\}\big\},\]
where the supremum is taken over all pairs $(g,J)$ that satisfy the conditions of boundedness (see page \pageref{bounded}). Here our convention is that $\sup\emptyset=0$. We define
\begin{eqnarray}\label{eq:Ham c M om}&\Ham_c(M,\om):=\big\{\phi_H^1\,\big|\,H\in C^\infty_c([0,1]\x M,\R)\big\},&\\
\label{eq:Vert om c}&\Vert\cdot\Vert_\om^c:\Ham_c(M,\om)\to[0,\infty),&\\
\nn&\Vert\phi\Vert_\om^c:=\inf\big\{\Vert H\Vert\,\big|\,H\in C^\infty_c([0,1]\x M,\R):\,\phi_H^1=\phi\big\},&
\end{eqnarray}
where $C^\infty_c([0,1]\x M,\R)$ denotes the space of all $H\in C^\infty([0,1]\x M,\R)$ with compact support. 
\begin{thm}[\cite{Ch}, Main Result]\label{thm:phi L} Let $(M,\om)$ be a symplectic manifold, $L\subseteq M$ a closed Lagrangian submanifold, and $\phi\in\Ham_c(M,\om)$. If 
\begin{equation}\label{eq:d c}\Vert\phi\Vert_\om^c<A_b(L)\end{equation}
then $\phi(L)\cap L\neq\emptyset$.
\end{thm}
\begin{rem}\label{rmk:phi L} The statement of Theorem \ref{thm:phi L} remains true if $\phi$ lies in the bigger group $\Ham(M,\om)$ and the condition (\ref{eq:d c}) is replaced by the weaker condition
\begin{equation}\label{eq:d}\Vert\phi\Vert_\om<A_b(L)
\end{equation}
This follows from Theorem \ref{thm:phi L} and Lemma \ref{le:phi psi} below. $\Box$
\end{rem}
\noindent{\bf Remark.} In Chekanov's Main Result it is assumed that $(M,\om)$ is bounded. This is unnecessary, since in the unbounded case we have $A_b(L)=0$, and hence the statement is void. $\Box$

\noindent{\bf Remark.} The definition of (geometric) boundedness in Y.~Chekanov's article is slightly stronger, and the number $A_b(L)$ in the hypothesis of the theorem is replaced by a corresponding quantity. However, the proof of the main result in that article goes through with these minor modifications. $\Box$

The proof of Theorem \ref{thm:N phi N} also relies on the following construction. Let $(M,\om)$ be a symplectic manifold, and $N\sub M$ a coisotropic submanifold. For $x\in N$ we denote by $N^\om_x$ the isotropic leaf of $N$ through $x$. Furthermore, we denote by $N_\om$ the set of isotropic leaves of $N$, and by $\pi_N:N\to N_\om$ the canonical projection. 

Assume that $N$ is regular. Then there exists a unique manifold structure on $N_\om$ such that $\pi_N$ is a smooth submersion. (This follows for example from \cite[Lemma 15, p.~20]{ZiLeafwise}.) Furthermore, there exists a unique symplectic structure $\om_N$ on $N_\om$ such that $\pi_N^*\om_N=\om|_N$. We define 
\begin{eqnarray}\label{eq:wt M wt om}
&\wt M:=M\x N_\om,\quad \wt\om:=\om\oplus(-\om_N),&\\
\label{eq:iota N wt N} &\iota_N:N\to \wt M,\quad \iota_N(x):=(x,N^\om_x),\quad \wt N:=\iota_N(N).&
\end{eqnarray} 
By a straight-forward argument the set $\wt N$ is a Lagrangian submanifold of $\wt M$. The next result is a crucial ingredient in the proof of Theorem \ref{thm:N phi N}. It is proved on page \pageref{proof:prop:A Cross wt M}. 
\begin{prop}\label{prop:A Cross wt M} If $N$ is closed and regular then 
\begin{equation}\label{eq:A Cross wt M} A_{\Cross}(M,\om,N)\leq A_{\Cross}(\wt M,\wt\om,\wt N).
\end{equation}
\end{prop}
We also need the following. Let $(M,\om)$ be a symplectic manifold and $L\sub M$ a Lagrangian submanifold. 
\begin{prop}\label{prop:A Cross b} If $M$ is connected and $L\neq\emptyset$ then we have
\[A_\Cross(M,\om,L)\leq A_b(M,\om,L).\]
\end{prop}
This result will be proved on page \pageref{proof:prop:A Cross b}. We are now ready for the proof of the key result. 

\begin{proof}[of Theorem \ref{thm:N phi N}]\setcounter{Claim}{0}\label{proof:thm:N phi N} Let $M,\om,N,\phi$ be as in the hypothesis, such that inequality (\ref{eq:d phi id}) is satisfied. Without loss of generality, we may assume that $M$ is connected. Consider the symplectomorphism
\[\hhat\phi:\wt M\to\wt M,\quad(x,x')\mapsto(\phi(x),x').\]
\begin{Claim}\label{claim:eq:d} We have
\begin{equation}\label{eq:Vert hhat phi}\Vert\hhat\phi\Vert_{\wt\om}<A_b(\wt M,\wt\om,\wt N).\end{equation}
\end{Claim}
\begin{proof}[of Claim \ref{claim:eq:d}] By a straight-forward argument we have that 
\begin{equation}\label{eq:d hhat phi id}\Vert\hhat\phi\Vert_{\wt\om}\leq\Vert\phi\Vert_\om.
\end{equation}
Since by hypothesis $N$ is regular and closed, we may apply Proposition \ref{prop:A Cross wt M}. It follows that inequality (\ref{eq:A Cross wt M}) holds. 

Furthermore, since by assumption $M$ and $N$ are connected, the manifold $\wt M$ is connected. Therefore, we may apply Proposition \ref{prop:A Cross b}, to conclude that
\begin{equation}\nn A_{\Cross}(\wt M,\wt\om,\wt N)\leq A_b(\wt M,\wt\om,\wt N).\end{equation}
Combining this with inequalities (\ref{eq:d hhat phi id},\ref{eq:d phi id},\ref{eq:A Cross wt M}), it follows that the inequality (\ref{eq:Vert hhat phi}) holds. This proves Claim \ref{claim:eq:d}.
\end{proof}
Since $N$ is closed, the manifold $\wt N$ is, as well. It follows that all hypotheses of Theorem \ref{thm:phi L} are satisfied, with $M,\om,\phi$ replaced by $\wt M,\wt\om,\hhat\phi$, and $L:=\wt N$, except for (\ref{eq:d c}). Furthermore, by Claim \ref{claim:eq:d}, the inequality (\ref{eq:d}) is satisfied. Therefore, using Remark \ref{rmk:phi L}, it follows that 
\begin{equation}\label{eq:hhat phi wt N}\hhat\phi(\wt N)\cap\wt N\neq\emptyset.\end{equation} 
We denote by $\pr:\wt M\to M$ the projection onto the first factor. Then we have
\begin{equation}\nn\pr\big(\hhat\phi(\wt N)\cap\wt N\big)\sub\phi(N)\cap N.\end{equation}
Combining this with (\ref{eq:hhat phi wt N}), the statement (\ref{eq:N phi N}) follows. This proves Theorem  \ref{thm:N phi N}.
\end{proof}
Next we will prove Proposition \ref{prop:A Cross wt M}. We will use the following construction. Let $(M,\om)$ and $(\hhat M,\hhat\om)$ be symplectic manifolds, $N\sub M$ and $\hhat N\sub \hhat M$ coisotropic submanifolds, and $\phi:\hhat M\to M$ a symplectomorphism satisfying $\phi(\hhat N)=N$. We define 
\begin{equation}\label{eq:phi'}\phi':\hhat N_{\hhat\om}\to N_\om,\quad\phi'(\hhat N^{\hhat\om}_{\hhat x}):=N^\om_{\phi(\hhat x)}.\end{equation} 
This map is well-defined. We also define 
\begin{equation}\label{eq:wt phi}\wt\phi:=\phi\x\phi':\wt{\hhat M}=\hhat M\x \hhat N_{\hhat\om}\to\wt M=M\x N_\om.
\end{equation}
\begin{rem}\label{rem:phi'} Assume that one of the manifolds $\hhat N$ or $N$ is regular. Then the other one is, as well, and $\phi'$ and hence $\wt\phi$ are symplectomorphisms. This follows from a straight-forward argument. $\Box$
\end{rem}
The proof of Proposition \ref{prop:A Cross wt M} also uses the following.
\begin{lem}\label{le:A wt A} Let $(M,\om)$ be a symplectic manifold and $N\sub M$ a closed, regular coisotropic submanifold. Then we have
\[A(M,\om,N)=A(\wt M,\wt\om,\wt N).\]
\end{lem}
\begin{proof}[of Lemma \ref{le:A wt A}]\setcounter{Claim}{0} This is Lemma 10 (Key Lemma) in \cite{ZiLeafwise}.
\end{proof}
\begin{proof}[of Proposition \ref{prop:A Cross wt M}]\setcounter{Claim}{0}\label{proof:prop:A Cross wt M} Assume that $(M_i,\om_i,N_i)_{i=1,\ldots,k}$ is a bounded splitting of $(M,\om,N)$. We define $(\wt{M_i},\wt{N_i})_{i=1,\ldots,k}$ as in (\ref{eq:wt M wt om},\ref{eq:iota N wt N}) with $M$ replaced by $M_i$ etc. By the definition of a bounded splitting of $(M,\om,N)$, there exists a symplectomorphism $\phi:\Cross_iM_i\to M$ such that $\phi(\Cross_iN_i)=N$. We define $\phi'$ and $\wt\phi$ as in (\ref{eq:phi'},\ref{eq:wt phi}), with 
\[(\hhat M,\hhat\om,\hhat N):=\big(\Cross_iM_i,\oplus_i\om_i,\Cross_iN_i\big).\] 
By hypothesis $N$ is regular. Hence by Remark \ref{rem:phi'}, the product $\Cross_iN_i$ is regular. Applying Lemma \ref{le:M 1 M 2} below, it follows that $N_i$ is regular, for every $i$. Since by hypothesis, $N$ is closed, Remark \ref{rem:X Y} below implies that $N_i$ is closed. We define the symplectic form $\wt\om_i$ on $\wt M_i$ as in (\ref{eq:iota N wt N}) with $M$ replaced by $M_i$ etc. 
\begin{Claim}\label{claim:splitting} The tuple $\big(\wt M_i,\wt\om_i,\wt N_i\big)_i$ is a bounded splitting of $(\wt M,\wt\om,\wt N)$. 
\end{Claim}
\begin{proof}[of Claim \ref{claim:splitting}] Let $i\in\{1,\ldots,k\}$. We show that {\bf $(\wt M_i,\wt\om_i)$ is bounded}: Since $N_i$ is closed, it follows that $(N_i)_{\om_i}$ is closed. Furthermore, by assumption, $(M_i,\om_i)$ is bounded. It follows that $(\wt M_i,\wt\om_i)$ is bounded.

We now show that {\bf there exists a map $\wt f:\Cross_i\wt M_i\to\wt M$} as in the definition of a bounded splitting: By Lemma \ref{le:products} below the identity map on $\Cross_iN_i$ descends to a symplectomorphism 
\[\psi':\Cross_i(N_i)_{\om_i}\to(\Cross_iN_i)_{\oplus_i\om_i}.\]
We denote by 
\[\wt\psi:\Cross_i\wt M_i=\Cross_i(M_i\x(N_i)_{\om_i})\to\Cross_iM_i\x(\Cross_iN_i)_{\oplus_i\om_i}\]
the map induced by $\psi'$, and define
\[\wt f:=\wt\phi\circ\wt\psi:\Cross_i\wt M_i\to\wt M.\] 
By Remark \ref{rem:phi'} the map $\wt\phi$ is a symplectomorphism. Since $\psi'$ is a symplectomorphism, the same holds for $\wt\psi$, and hence for $\wt f$. Furthermore, we have $\wt f(\Cross_i\wt N_i)=\wt N$. Hence the map $\wt f$ satisfies the conditions in the definition of a bounded splitting. This proves Claim \ref{claim:splitting}.
\end{proof}
Let $i=1,\ldots,k$. Since $N_i$ is closed, we may apply Lemma \ref{le:A wt A}, to conclude that $A\big(\wt{M_i},\wt{\om_i},\wt{N_i}\big)=A(M_i,\om_i,N_i)$. Combining this with Claim \ref{claim:splitting}, the inequality (\ref{eq:A Cross wt M}) follows. This proves Proposition \ref{prop:A Cross wt M}.
\end{proof}
For the proof of Proposition \ref{prop:A Cross b} we need the following. Recall the definition (\ref{eq:A Si M J X}). 
\begin{lem}\label{le:M i om i X i} Let $\Si$ be a Riemann surface, $k\in\N$, and for $i=1,\ldots,k$ let $(M_i,\om_i)$ be a symplectic manifold, $X_i\sub M_i$ a subset, and $J_i$ an $\om_i$-tame almost complex structure on $M_i$. Then 
\begin{equation}\label{eq:min A}\min_{i=1,\ldots,k}A\big(\Si,M_i,\om_i,J_i,X_i\big)\leq A\big(\Si,\Cross_iM_i,\oplus_i\om_i,\oplus_iJ_i,\Cross_iX_i\big).\end{equation}
\end{lem}
In the proof of this lemma we will use the following. 
\begin{rem}\label{rem:pos} Let $(M,\om)$ be a symplectic manifold, $J$ an $\om$-tame almost complex structure, $\Si$ a Riemann surface, and $u:\Si\to M$ a $J$-holomorphic map. Then $\int_\Si u^*\om\geq0$. This follows from the fact that $\int_\Si u^*\om$ is the Dirichlet energy of $u$. $\Box$
\end{rem}
\begin{proof}[of Lemma \ref{le:M i om i X i}]\setcounter{Claim}{0} Assume that $u:=(u_1,\ldots,u_k):\Si\to M:=\Cross_iM_i$ is a $J:=\oplus_iJ_i$-holomorphic map satisfying $u(\dd\Si)\sub X:=\Cross_iX_i$ and $E:=\int_\Si u^*\om>0$. Let $i=1,\ldots,k$. We denote $E_i:=\int_\Si u_i^*\om_i$. Since by assumption $J_i$ is $\om_i$-tame, and $u_i$ is $J_i$-holomorphic, by Remark \ref{rem:pos} we have $E_i\geq0$. Combining this with the fact $E=\sum_iE_i$, it follows that 
\begin{equation}\label{eq:E i E}E_i\leq E,\quad\forall i=1,\ldots,k.\end{equation} 
Since $E>0$, there exists $i_0\in\{1,\ldots,k\}$ such that $E_{i_0}>0$. Combining this with inequality (\ref{eq:E i E}) and using the fact $u_i(\dd\Si)\sub X_i$, the inequality (\ref{eq:min A}) follows. This proves Lemma \ref{le:M i om i X i}.
\end{proof}
In the proof of Proposition \ref{prop:A Cross b} we will use the following remark. We define the bounded minimal action as in (\ref{eq:A b}).
\begin{rem}\label{rmk:A b} Let $(M,\om)$ be a symplectic manifold, $N\sub M$ a coisotropic submanifold, $M'$ a smooth manifold, and $\phi:M'\to M$ a diffeomorphism. Then we have
\[A_b\big(M',\phi^*\om,\phi^{-1}(N)\big)=A_b(M,\om,N).\]
\end{rem}
\begin{proof}[of Proposition \ref{prop:A Cross b}]\setcounter{Claim}{0}\label{proof:prop:A Cross b} Let $(M_i,\om_i,L_i)_{i=1,\ldots,k}$ be a bounded splitting of $(M,\om,L)$. The statement of the proposition is a consequence of the following claim:
\begin{equation}\label{eq:min i a i}\min_iA(M_i,\om_i,L_i)\leq A_b(M,\om,L).\end{equation}
To see that this inequality holds, we choose a map $\phi$ as in the definition of a bounded splitting. By Remark \ref{rmk:A b} we may assume without loss of generality that $(M,\om)=(\Cross_iM_i,\oplus_i\om_i)$ and $\phi=\id$. For $i=1,\ldots,k$ we choose a pair $(g_i,J_i)$ as in the definition of boundedness of $(M_i,\om_i)$, and we define
\[a_i:=\min\big\{A\big(\D,M_i,\om_i,J_i,L_i\big),A\big(S^2,M_i,\om_i,J_i,\emptyset\big)\big\}.\]
Let $i\in\{1,\ldots,k\}$. The submanifold $L_i\sub M_i$ is Lagrangian. Furthermore, since by hypothesis $M$ is connected and non-empty, the manifold $M_i$ is connected. Hence, using Lemma \ref{le:S M om N S 2} below, we have 
\begin{equation}\label{eq:a i} A(M_i,\om_i,L_i)\leq a_i.
\end{equation}
We define $(g,J):=\big(\oplus_ig_i,\oplus_iJ_i\big)$. This pair satisfies the conditions of boundedness for $(M,\om)$. Therefore, we have
\begin{equation}\label{eq:A D S 2}a:=\min\big\{A\big(\D,M,\om,J,L\big),A\big(S^2,M,\om,J,\emptyset\big)\big\}\leq A_b(M,\om,L).\end{equation}
It follows from the definition of boundedness that $J_i$ is $\om_i$-tame, for every $i$. Therefore, we may apply Lemma \ref{le:M i om i X i}. It follows that $\min_ia_i\leq a$. Combining this with (\ref{eq:A D S 2},\ref{eq:a i}), inequality (\ref{eq:min i a i}) follows. This completes the proof of Proposition \ref{prop:A Cross b}. \end{proof}
\subsection{Proofs of Theorem \ref{thm:A coiso cap} (Regular coisotropic capacity) and Corollary \ref{cor:e w} (Energy-Gromov-width inequality)}\label{proof:A coiso cap,cor:e w}
We will first prove Theorem \ref{thm:A coiso cap}, since it is used in the proof of Corollary \ref{cor:e w}. For the proofs of both results we need the following lemma.
\begin{lem}\label{le:A coiso Cross d} Let $(M,\om)$ be a bounded and aspherical symplectic manifold of dimension $2n$, $U\sub M$ an open subset, and $d\in\{n,\ldots,2n\}$. Then we have
\begin{equation}\label{eq:A coiso Cross d}A^\coiso_\Cross(M,\om,U)\geq A_\coiso^d(U,\om|_U).\end{equation}
\end{lem}
\begin{proof}[of Lemma \ref{le:A coiso Cross d}] Let $N\sub U$ be a regular closed coisotropic submanifold of dimension $d$ such that condition (\ref{eq:contr}) is satisfied. Since by hypothesis $(M,\om)$ is aspherical, Lemma \ref{lemactionspec} below implies that 
\begin{equation}\label{eq:A N A U}A(M,\om,N)\geq A(U,\om|_U,N).\end{equation}
Since by hypothesis $(M,\om)$ is bounded, we have $A^\coiso_\Cross(M,\om,N)\geq A(M,\om,N)$. Combining this with inequality (\ref{eq:A N A U}), inequality (\ref{eq:A coiso Cross d}) follows. This proves Lemma \ref{le:A coiso Cross d}.
\end{proof}

The proof of statement (\ref{thm:A coiso cap:A coiso d B Z}) of Theorem \ref{thm:A coiso cap} involves a certain product of rescaled Stiefel manifolds. These manifolds are given as follows. Let $k,n\in\N$ be such that $k\leq n$, and $a>0$. We define the \emph{Stiefel manifold of symplectic area $a$} to be
\[V(k,n,a):=\big\{\Theta\in\C^{k\x n}\,\big|\,\Theta\Theta^*=\frac a\pi\one_k\big\}.\]
The proofs of statements (\ref{thm:A coiso cap:cap},\ref{thm:A coiso cap:asph})) involve the spherical action (or area) spectrum of a symplectic manifold $(M,\om)$. It is given by 
\begin{equation}\label{eq:S M om} S(M,\om):=\inf\big(\big\{\int_{S^2}u^*\om\,\big|\,u\in C^\infty(S^2,M)\big\}\cap(0,\infty)\big)\in[0,\infty].
\end{equation} 

\begin{proof}[of Theorem \ref{thm:A coiso cap}] \setcounter{Claim}{0} We start by proving {\bf statement (\ref{thm:A coiso cap:A coiso d B Z})}. To see that (\ref{eq:A coiso d B}) holds, let $\ell\in\N$ and $k_i,n_i$, $i=1,\ldots,\ell$ be as in the definition of $\keq(n,d)$ (see (\ref{eq:k})). We define $a_0:=\pi/\sum_{i=1}^\ell k_i$. Let $a\in(0,a_0)$ be a number. We define $N:=\Cross_{i=1}^\ell V(k_i,n_i,a)$. Then $N$ is a closed regular coisotropic submanifold of $B^{2n}\sub\R^{2n}=\Cross_i\C^{k_i\x n_i}$ of dimension $d$. Since $\R^{2n}$ is simply connected, every loop in an isotropic leaf of $N$ is contractible in $\R^{2n}$. 

Hence $N$ satisfies the conditions in the definition of $A_\coiso^d(B^{2n},\om_0)$. Furthermore, by Lemma \ref{lemactionspec}, Remark \ref{rmk:S N N'}, and Proposition \ref{minimalareaStiefel} below we have 
\[A^d_\coiso(B^{2n},\om_0)\geq A(B^{2n},\om_0,N)=A(\R^{2n},\om_0,N)=a.\] 
Since this holds for arbitrary $a\in(0,a_0)$, it follows that $A^d_\coiso(B^{2n},\om_0)\geq a_0$. Inequality (\ref{eq:A coiso d B}) follows. 

To see that (\ref{eq:A coiso d Z}) holds, note that by Lemma \ref{le:A coiso Cross d}, we have
\[A_\coiso^d(Z^{2n},\om_0)\leq A^\coiso_\Cross(\R^{2n},\om_0,Z^{2n}).\]
Combining this with Corollary \ref{cor:e A Cross} and Remark \ref{rem:e Z 2n a} below, the inequality (\ref{eq:A coiso d Z}) follows. 

We prove {\bf statement (\ref{thm:A coiso cap:cap})}. Let $n\in\N$ and $d\in\{n,\ldots,2n-1\}$. To prove {\bf{monotonicity}} of the restriction of $A_\coiso^d$, let $(M,\om)$ and $(M',\om')$ be aspherical symplectic manifolds of dimension $2n$, and $\phi:M'\to M$ a symplectic embedding. It follows from Lemma \ref{lemactionspec} below and asphericity of $(M,\om)$, that for every regular closed coisotropic submanifold $N'\sub M'$ satisfying (\ref{eq:contr}), we have
\[A\big(M,\om,\phi(N')\big)\geq A(M',\om',N').\]
It follows that
\[A_\coiso^d(M',\om')\leq A_\coiso^d(M,\om).\]
This proves (monotonicity).

{\bf{Conformality}} follows immediately from the definitions. \\
{\bf{Non-triviality}} follows from the inequalities (\ref{eq:A coiso d B},\ref{eq:A coiso d Z}). Furthermore, inequality (\ref{eq:2n d}) in Proposition \ref{prop:k n d} implies that $\keq(n,2n-1)=1$. Hence it follows from inequalities (\ref{eq:A coiso d B},\ref{eq:A coiso d Z}) that $A_\coiso^{2n-1}$ is normalized. This proves statement (\ref{thm:A coiso cap:cap}).

To prove {\bf statement (\ref{thm:A coiso cap:asph})}, let $N\sub M$ be a closed regular coisotropic submanifold of dimension $d$. We define $\wt N:=N\x M'$. This is a closed and regular coisotropic submanifold of $\wt M:=M\x M'$, of dimension $d+2n'$. By Lemma \ref{le:S M} below and asphericity of $(M',\om')$ we have 
\[S(M',\om',M')\sub S(M',\om')=\{0\}.\] 
Hence Remark \ref{rmk:S N N'} below implies that $A(\wt N)=A(N)$. Furthermore, if $N$ satisfies (\ref{eq:contr}) then the same holds for $\wt N$. Hence the inequality (\ref{eq:A coiso d 2n'}) follows. This proves statement (\ref{thm:A coiso cap:asph}) and completes the proof of Theorem \ref{thm:A coiso cap}.
\end{proof}
We are now ready to prove Corollary \ref{cor:e w}. 

\begin{proof}[of Corollary \ref{cor:e w}]\setcounter{Claim}{0}\label{proof:cor:e w} By Corollary \ref{cor:e A Cross} we have
\begin{equation}\label{eq:e U M' A} e\big(U\x M',M\x M',\om\oplus\om'\big)\geq A^\coiso_\Cross\big(M\x M',\om\oplus\om',U\x M'\big).
\end{equation}
We denote $2n:=\dim M$ and $2n':=\dim M'$. Since by hypothesis $(M,\om)$ and $(M',\om')$ are bounded and aspherical, the same holds for their product. Hence applying Lemma \ref{le:A coiso Cross d}, we obtain
\begin{equation}\label{eq:A coiso Cross M M'}A^\coiso_\Cross\big(M\x M',\om\oplus\om',U\x M'\big)\geq A_\coiso^{2(n+n')-1}\big(U\x M',\om|_U\oplus\om'\big).
\end{equation}
Using closedness and asphericity of $(M',\om')$, Theorem \ref{thm:A coiso cap}(\ref{thm:A coiso cap:asph}) implies that
\begin{equation}\label{eq:A coiso U M'}A_\coiso^{2(n+n')-1}\big(U\x M',\om|_U\oplus\om'\big)\geq A_\coiso^{2n-1}(U,\om|_U).
\end{equation}
Using asphericity of $U$, Theorem \ref{thm:A coiso cap}(\ref{thm:A coiso cap:cap}) implies that
\[A_\coiso^{2n-1}(U,\om|_U)\geq w(U,\om|_U).\]
Combining this with inequalities (\ref{eq:e U M' A},\ref{eq:A coiso Cross M M'},\ref{eq:A coiso U M'}), the inequality (\ref{eq:e w}) follows. This proves Corollary \ref{cor:e w}.
\end{proof}
\subsection{Proof of Corollary \ref{cor:skinny} (Skinny non-squeezing)} \label{proof:cor:skinny,q:M U, q:M U a} 
\begin{proof}[of Corollary \ref{cor:skinny}]\setcounter{Claim}{0}\label{proof:cor:skinny} We denote 
\[\wt M:=\R^{2n}\x M,\quad\wt\om:=\om_0\oplus\om.\]
Let $\phi:U\x M\to \wt M$ be a symplectic embedding. It suffices to prove the following. Assume that $a_0>0$ is such that 
\begin{equation}\label{eq:phi U M}\phi(U\x M)\sub Z^{2n}(a_0)\x M.\end{equation}
Then we have 
\begin{equation}\label{eq:a 0 pi}a_0>\pi.\end{equation}
To see that this inequality holds, we define $N:=\Cross_{i=1}^kS^{2n_i-1}\x M$ and $N':=\phi(N)$. 
\begin{Claim}\label{claim:a N'}Let $a>0$ be such that $N'\sub Z^{2n}(a)\x M$. Then we have $a\geq\pi$.
\end{Claim}
\begin{proof}[of Claim \ref{claim:a N'}] The inclusion $N'\sub Z^{2n}(a)\x M$ and Remarks \ref{rem:e M M'} and \ref{rem:e Z 2n a} below imply that
\[e(N',\wt M,\wt\om)\leq a.\]
Therefore, the inequality $a\geq\pi$ is a consequence of the following claim. 
\begin{Claim}\label{claim:e N' wt M} We have
\begin{equation}\label{eq:e N' wt M}e(N',\wt M,\wt\om)\geq\pi.\end{equation} 
\end{Claim}
\begin{pf}[Proof (of Claim \ref{claim:e N' wt M})] $N$ is a closed and regular coisotropic submanifold of $U\x M$, and hence $N'$ is a closed and regular coisotropic submanifold of $\wt M$. We define $n:=\sum_{i=1}^kn_i$. By Remark \ref{rmk:S N N'} below we have
\begin{eqnarray}\nn S\big(U\x M,\wt\om,N\big)&\sub&S\big(\R^{2n}\x M,\wt\om,N\big)\\
\label{eq:S U M}&=&\Big(\sum_iS\big(\R^{2n_i},\om_0,S^{2n_i-1}\big)\Big)+S(M,\om,M).\end{eqnarray}
Proposition \ref{minimalareaStiefel} below implies that $S\big(\R^{2n_i},\om_0,S^{2n_i-1}\big)=\pi\Z$. Furthermore, by Lemma \ref{le:S M om N S 2}(\ref{le:S M om N S 2:M}) below and the hypothesis (\ref{eq:int u om}), we have $S(M,\om,M)=S(M,\om)\sub\pi\Z$. Combining this with (\ref{eq:S U M}), it follows that 
\begin{equation}\label{eq:S U M pi Z}S\big(U\x M,\wt\om,N\big)\sub\pi\Z.\end{equation} 
By Lemma \ref{le:leaf prod} below the isotropic leaves of $N$ are the products of the isotropic leaves of $\Cross_iS^{2n_i-1}$ and $M$ (viewed as a coisotropic submanifold of itself). The latter are single points. Furthermore, the hypothesis $n_i\geq2$, for every $i$, implies that $\Cross_iS^{2n_i-1}$ is simply-connected. It follows that every loop in an isotropic leaf of $N$ is contractible in $N$, and hence in $U$. Hence we may apply Lemma \ref{lemactionspec} below, and conclude that 
\begin{equation}\label{eq:S R 2n M}S(\wt M,\wt\om,N')\sub S\big(U\x M,\wt\om,N\big)+S(\wt M,\wt\om).\end{equation}
The hypothesis (\ref{eq:int u om}) implies that $S(\wt M,\wt\om)\sub\pi\Z$. Combining this with (\ref{eq:S U M pi Z},\ref{eq:S R 2n M}), it follows that $S(\wt M,\wt\om,N')\sub\pi\Z$, and therefore,
\begin{equation}\label{eq:A wt M}A\big(\wt M,\wt\om,N'\big)\geq\pi.\end{equation}
Since by hypothesis $M$ is closed, the symplectic manifold $(M,\om)$ is bounded. Hence the same holds for $(\wt M,\wt\om)$. It follows that 
\[A_\Cross(\wt M,\wt\om,N')\geq A\big(\wt M,\wt\om,N'\big).\]
Combining this with (\ref{eq:A wt M}) and applying Theorem \ref{thm:N phi N}, inequality (\ref{eq:e N' wt M}) follows. This proves Claim \ref{claim:e N' wt M} and hence Claim \ref{claim:a N'}.
\end{pf}
\end{proof}
Using the assumption (\ref{eq:phi U M}) and compactness of $N'$, there exists $a<a_0$ such that $N'\sub Z^{2n}(a)\x M$. By Claim \ref{claim:a N'}, it follows that $a\geq\pi$. Inequality (\ref{eq:a 0 pi}) follows. This proves Corollary \ref{cor:skinny}.
\end{proof}
\subsection{Proof of Theorem \ref{thm:emb d B} (Badly squeezable small sets)}\setcounter{Claim}{0}\label{subsec:proof:thm:emb d B}
For the proof of this theorem, we need the following results. Let $d\in[0,\infty)$. For every $n\in\N$ we abbreviate
\[\emb_d:=\emb^{Z^{2n},\om_0}_d.\]
Let $(M,\om)$ be a symplectic manifold.
\begin{prop}\label{prop:c M om d} For every $n\in\N$, we have 
\begin{equation}\label{eq:emb d M om 2n'}\emb_d(M,\om)\geq \emb_d(M\x\R^{2n},\om\oplus\om_0).\end{equation}
\end{prop}
We post-pone the proof of this result to page \pageref{proof:prop:c M om d}.
\begin{prop}\label{prop:ineqembA} If $d\geq2$ then we have
\begin{equation}\label{eq:emb d A}
\emb_d(M,\omega)\geq A_\coiso^d(M,\omega).
\end{equation}
\end{prop}
We post-pone the proof of this result to page \pageref{proof:prop:ineqembA}. For the proof of Theorem \ref{thm:emb d B} we also need the following. 
\begin{rem}\label{rem:scaling} Let $n\in\N$, $d\in[0,\infty)$, $r>0$, and $U\sub \R^{2n}$ an open subset. Then 
\[\emb_d(rU,\om_0)=r^2\emb_d(U,\om_0).\]
This follows from a straight-forward argument. $\Box$
\end{rem}
For $k,n\in\N$ satisfying $k\leq n$ we denote by 
\[V(k,n):=\big\{\Theta\in\C^{k\x n}\,\big|\,\Theta\Theta^*=\one_k\big\}\]
the Stiefel manifold of unitary $k$-frames in $\C^n$. 

\begin{proof}[of Theorem \ref{thm:emb d B}]\setcounter{Claim}{0} Let $n\in\{2,3\ldots\}$ and $d\in[n,\infty)$. By Remark \ref{rem:scaling} it suffices to prove that 
\begin{equation}\label{eq:emb d n}\emb_d(B^{2n}(a),\om_0)\geq\pi,\quad\forall a>\pi\kk(n,d).\end{equation}
To show that this condition holds, let $a>\pi\kk(n,d)$. We choose $\ell\in\N$ and $k_1,\ldots,k_\ell\in\N$ such that $\sum_{i=1}^\ell k_i=\kk(n,d)$. We also choose $n_1,\ldots,n_{\ell}$ satisfying (\ref{eq:n i k i},\ref{eq:k i n i d geq}). We define $n':=\sum_ik_in_i-n$. By the first inequality in (\ref{eq:k i n i d geq}) we have $n'\geq0$. Propositions \ref{prop:c M om d} and \ref{prop:ineqembA} imply that 
\begin{eqnarray}
\nn\emb_d(B^{2n}(a),\omega_0)&\geq&\emb_d(B^{2n}(a)\times\R^{2n'},\omega_0)\\
\label{eq:A coiso d B 2n a}&\geq& A_\coiso^d(B^{2n}(a)\times\R^{2n'},\omega_0).
\end{eqnarray}
By the inequalities (\ref{eq:n i k i}) the Stiefel manifolds $V(k_i,n_i)$ are well-defined. We define $N:=V(k_1,n_1)\x\ldots\x V(k_\ell,n_\ell)$. 
\begin{Claim}\label{claim:A coiso d} We have
\begin{equation}\label{eq:A coiso d B R}A_\coiso^d(B^{2n}(a)\times\R^{2n'},\omega_0)\geq A(\R^{2(n+n')},\omega_0,N).\end{equation}
\end{Claim}
\begin{proof}[of Claim \ref{claim:A coiso d}] Note that $N$ is a regular closed coisotropic submanifold of $B^{2n}(a)\x\R^{2n'}$ of dimension $\sum_ik_i(2n_i-k_i)$. Hence by the second inequality in (\ref{eq:k i n i d geq}) we have $\dim N\leq d$. Furthermore, since $B^{2n}(a)\x\R^{2n'}$ is contractible, condition (\ref{eq:contr}) is satisfied. Therefore, by definition, we have 
\[A_\coiso^d(B^{2n}(a)\times\R^{2n'},\omega_0)\geq A(B^{2n}(a)\times\R^{2n'},\omega_0,N).\]
The right hand side is bounded below by $A(\R^{2(n+n')},\omega_0,N)$. Inequality (\ref{eq:A coiso d B R}) follows. This proves Claim \ref{claim:A coiso d}.
\end{proof}
By Remark \ref{rmk:S N N'} and Proposition \ref{minimalareaStiefel} below we have 
\[A(\R^{2(n+n')},\omega_0,N)=\pi.\]
Combining this with (\ref{eq:A coiso d B 2n a}) and Claim \ref{claim:A coiso d}, inequality (\ref{eq:emb d n}) follows. This proves Theorem \ref{thm:emb d B}.
\end{proof}
\begin{proof}[of Proposition \ref{prop:c M om d}]\setcounter{Claim}{0}\label{proof:prop:c M om d} Let $\wt X\sub M\x\R^{2n}$ be a compact subset of Hausdorff dimension at most $d$. We denote $2m:=\dim M$, $\wt M:=M\x\R^{2n}$, $\wt\om:=\om\oplus\om_0$, by $\pr:\wt M\to M$ the projection onto the first component, and $X:=\pr(\wt X)$. Then $X$ is a compact subset of $M$. Furthermore, by standard results (cf.~\cite{Fed}), the Hausdorff dimension of $X$ does not exceed that of $\wt X$, and thus is at most $d$. Recall the definition (\ref{eq:emb M 0 om 0 M om X}). 
\begin{Claim}\label{claim:emb X X'} We have 
\begin{equation}\label{eq:c M om X}\emb^{Z^{2m},\om_0}(M,\om,X)\geq\emb^{Z^{2(m+n)},\om_0}(\wt M,\wt\om,\wt X).
\end{equation}
\end{Claim}
\begin{proof}[of Claim \ref{claim:emb X X'}] Let $a>0$. Assume that there exists a pair $(U,\phi)$, where $U\sub M$ is an open neighborhood of $X$ and $\phi:U\to Z^{2m}(a)$ is a symplectic embedding. We define $\wt U:=U\x\R^{2n}$ and $\wt\phi:=\phi\x\id_{\R^{2n}}$. Then $\wt X\sub\wt U$ and $\wt\phi$ is a symplectic embedding of $\wt U$ into $Z^{2(m+n)}(a)$. The inequality (\ref{eq:c M om X}) follows. This proves Claim \ref{claim:emb X X'}.
\end{proof}
Taking the supremum over all compact sets $\wt X\sub M\x\R^{2n}$ of Hausdorff dimension at most $d$, Claim \ref{claim:emb X X'} implies inequality (\ref{eq:emb d M om 2n'}). This completes the proof of Proposition \ref{prop:c M om d}.
\end{proof}
The idea of proof of Proposition \ref{prop:ineqembA} is the following. Let $N\sub M$ be a $d$-dimensional closed regular coisotropic submanifold satisfying (\ref{eq:contr}). We glue finitely many disks to $N$, in such a way that every loop in an isotropic fiber of $N$ is contractible in the resulting subset of $M$. This is possible because of (\ref{eq:contr}) and regularity and closedness of $N$. The statement of Proposition \ref{prop:ineqembA} will then be a consequence of Theorem \ref{thm:N phi N}, Lemma \ref{lemactionspec}, and Remark \ref{rem:e Z 2n a} below.
\begin{proof}[of Proposition \ref{prop:ineqembA}]\setcounter{Claim}{0}\label{proof:prop:ineqembA} Let $N\subseteq M$ be a non-empty closed regular coisotropic submanifold of dimension at most $d$, satisfying (\ref{eq:contr}). Inequality (\ref{eq:emb d A}) is a consequence of the following claim.
\begin{Claim}\label{claim:emb d A N} We have
\begin{equation}\label{eq:emb d A N}\emb_d(M,\om)\geq A(N).
\end{equation}
\end{Claim}
\begin{pf}[Proof (of Claim \ref{claim:emb d A N})] Without loss of generality, we may assume that $N$ is connected. We choose an isotropic leaf $F\sub N$ and a point $x_0\in F$. Regularity of $N$ implies that $F$ is a smooth submanifold of $N$. It is closed, since $N$ is closed. It follows that the fundamental group of $F$ with base point $x_0$ is finitely generated. Therefore, there exists a finite set $S$ of smooth loops $x:S^1\sub\C\to F$ satisfying $x(1)=x_0$, whose continuous homotopy classes with fixed base point generate $\pi_1(F,x_0)$. 

The assumption (\ref{eq:contr}) implies that for every $x\in S$ there exists a smooth map $u_x:\D\to M$ satisfying $u_x|_{S^1}=x$. We choose such a collection of maps $(u_x)_{x\in S}$ and define
\[X:=N\cup\bigcup_{x\in S}u_x(\D)\sub M.\]
This set is compact. Furthermore, a standard result (cf.~\cite[p.~176]{Fed}) implies that $u_x(\D)$ has Hausdorff dimension at most 2. Since by hypothesis $d\geq2$, it follows that $X$ has Hausdorff dimension at most $d$. We denote $2n:=\dim M$. Assume that $a>0$ is such that there exists a pair $(U,\phi)$, where $U\sub M$ is an open neighborhood of $X$ and $\phi:U\inj Z^{2n}(a)$ is a symplectic embedding. Using the fact $A(N)=A(M,\om,N)\leq A(U,\om|_U,N)$, inequality (\ref{eq:emb d A N}) is a consequence of the following claim. 
\begin{Claim}\label{claim:A a} We have 
\begin{equation}\label{eq:A a}a\geq A(U,\om|_U,N).
\end{equation}
\end{Claim}
\begin{pf}[Proof (of Claim \ref{claim:A a})] We choose a pair $(U,\phi)$ as above. 
\begin{Claim}\label{claim:loop N} Every continuous loop in an isotropic leaf of $N$ is contractible in $X$. 
\end{Claim}
In the proof of this claim we use the following notation. Let $X$ be a set and $x:S^1\sub\C\to X$ a map. We define $x^{-1}:S^1\to X$ by $x^{-1}(z):=x(\BAR z)$.

\begin{proof}[of Claim \ref{claim:loop N}] Let $x$ be such a loop. {\bf Assume first that $x(S^1)\sub F$}. It follows from our choice of the set $S$ that there exist $\ell\in\N_0$, $x_1,\ldots,x_\ell\in S$, and $\eps_1,\ldots,\eps_\ell\in\{1,-1\}$, such that $x$ is continuously homotopic inside $F$ to $x_1^{\eps_1}\cdots x_\ell^{\eps_\ell}$. Since $X$ contains the images $u_{x_i}(\D)$, for $i=1,\ldots,\ell$, it follows that $x$ is contractible in $X$. 

{\bf Consider now the general situation}. Since $N$ is path-connected, the same holds for the set of isotropic leaves $N_\om$. Hence there exists a path $\BAR{y}\in C([0,1],N_\om)$ such that $\BAR{y}(0)$ is the leaf through $x_0$, and $\BAR{y}(1)$ is the leaf containing $x(S^1)$. 

By regularity of $N$ there exists a unique manifold structure on $N_\om$, such that the canonical projection $\pi_N:N\to N_\om$ is a smooth submersion. (See Lemma 15 in \cite{ZiLeafwise}.) Since $N$ is closed, C.~Ehresmann's Theorem implies that $\pi_N$ is a smooth (locally trivial) fiber bundle. It follows that $\pi_N$ has the continuous homotopy lifting property. Hence there exists $u\in C([0,1]\x S^1,N)$ such that $\pr_N\circ u(t,z)=\bar y(t)$, for every $t\in[0,1]$ and $z\in S^1$, and $u(1,\cdot)=x$. By what we already proved, the loop $u(0,\cdot)$ is contractible in $X$. It follows that the same holds for $u(1,\cdot)=x$. This proves Claim \ref{claim:loop N}.
\end{proof}
Using Claim \ref{claim:loop N} and asphericity of $(\R^{2n},\om_0)$, Lemma \ref{lemactionspec} implies that 
\begin{equation}\label{eq:A phi N}A(U,\om|_U,N)\leq A(\R^{2n},\om_0,\phi(N)).
\end{equation}
Furthermore, the coisotropic submanifold $\phi(N)\sub\R^{2n}$ is non-empty, closed, and regular. Hence Theorem \ref{thm:N phi N} implies that 
\begin{equation}\label{eq:A phi N e}A(\R^{2n},\om_0,\phi(N))\leq e\big(\phi(N),\R^{2n},\om_0\big).\end{equation}
By Remark \ref{rem:e Z 2n a}, we have
\[e\big(\phi(N),\R^{2n},\om_0\big)\leq e\big(Z^{2n}(a),\om_0,\R^{2n},\om_0\big)\leq a.\]
Combining this with (\ref{eq:A phi N},\ref{eq:A phi N e}), inequality (\ref{eq:A a}) follows. This proves Claim \ref{claim:A a} and hence Claim \ref{claim:emb d A N}, and concludes the proof of Proposition \ref{prop:ineqembA}. 
\end{pf}\end{pf}\end{proof}
\subsection{Proof of Corollary \ref{cor:exotic} (Stably exotic form)}\label{proof:cor:exotic}
We need the following results, in which $(M,\om)$ and $(M',\om')$ are symplectic manifolds. 
\begin{cor}\label{cor}\label{cor:infinite} Assume that $(M,\om)$ is (geometrically) bounded and every compact subset of $M$ is Hamiltonianly displaceable. Then $(M,\om)$ is \emph{not} (coisotropically) infinite.
\end{cor}
\begin{proof}[of Corollary \ref{cor:infinite}]\setcounter{Claim}{0} This is a direct consequence of Theorem \ref{thm:N phi N}.
\end{proof}
We call $(M,\om)$ \emph{strongly (coisotropically) infinite} iff there exists a non-empty regular closed coisotropic submanifold $N\sub M$ such that $A(N)=\infty$, and every continuous loop in an isotropic leaf of $N$ is contractible in $M$.
\begin{prop}\label{prop:infinite} The following statements hold.
\begin{enui} 
\item\label{prop:infinite:asph}(Aspherical manifold) If $(M,\om)$ is closed and aspherical then it is strongly infinite. 
\item\label{prop:infinite:prod}(Product) The product of two infinite/ strongly infinite symplectic manifolds is infinite/ strongly infinite.
\item\label{prop:infinite:emb}(Embedding) Assume that $(M,\om)$ is aspherical and $(M',\om')$ is strongly infinite and embeds into $(M,\om)$. Then $(M,\om)$ is strongly infinite.
\end{enui}
\end{prop} 
\begin{proof}[of Proposition \ref{prop:infinite}]\setcounter{Claim}{0} {\bf Statement (\ref{prop:infinite:asph})} follows from Lemma \ref{le:S M om N S 2}(\ref{le:S M om N S 2:M}) below, using that $N:=M$ is a regular coisotropic submanifold of itself and the fact that the isotropic leaves of $M$ are single points. {\bf Statement (\ref{prop:infinite:prod})} follows from Remark \ref{rmk:S N N'} and Lemma \ref{le:leaf prod} below. {\bf Statement (\ref{prop:infinite:emb})} follows from Lemma \ref{lemactionspec} below.
\end{proof}
\begin{proof}[of Corollary \ref{cor:exotic}]  \setcounter{Claim}{0}
Assume that $\om$ is a (coisotropically) infinite form on $\R^{2n}$. Since $\R^{2n}$ is simply connected, it follows that $\om$ is strongly infinite. Let $X,\si,\si',\phi$ be as in the definition of stable exoticness. We show that condition (\ref{eq:phi si om 0}) holds. Consider first the {\bf case} in which $\si'$ is \emph{not} aspherical. Then (\ref{eq:phi si om 0}) holds, since $\si\oplus\om_0$ is aspherical. 

Consider now the {\bf case} in which $\si'$ \emph{is} aspherical. Then it follows from Proposition \ref{prop:infinite}(\ref{prop:infinite:asph},\ref{prop:infinite:prod}) that $(M,\Om'):=\big(X\x\R^{2n},\si'\oplus\om\big)$ is strongly infinite. We define $\Om:=\si\oplus\om_0$. The symplectic manifold $(M,\Om)$ is bounded, since it is the product of two bounded symplectic manifolds. Furthermore, every compact subset of $M$ is displaceable in an $\Om$-Hamiltonian way, since $(\R^{2n},\om_0)$ has this property. Therefore, by Corollary \ref{cor:infinite}, $(M,\Om)$ is not infinite. Moreover, since, by assumption, $\si$ is aspherical, $\Om$ is aspherical. 

Combining these facts and using that $(M,\Om')$ is strongly infinite, it follows from Proposition \ref{prop:infinite}(\ref{prop:infinite:emb}) that $(M,\Om')$ does not embed into $(M,\Om)$. Therefore, Condition (\ref{eq:phi si om 0}) holds. This completes the proof of Corollary \ref{cor:exotic}.
\end{proof}
\appendix
\section{Auxiliary results}\label{sec:aux}
\subsection{(Pre-)symplectic geometry}\label{subsec:pre sympl}
The next result is used in the proofs of Lemmas \ref{le:leaf prod} and \ref{le:M 1 M 2} below. 

Let $V$ be a finite dimensional vector space and $\om$ a skew-symmetric 2-form on $V$. We define 
\[V^\om:=\big\{v\in V\,\big|\,\om(v,w)=0,\,\forall w\in V\big\},\quad\corank\om:=\dim V^\om.\] 
By a \emph{presymplectic structure} on a manifold $M$ we mean a closed two-form $\om$ on $M$ such that $\corank\om_x$ does not depend on $x\in M$. Note that if $(M,\om)$ is a symplectic manifold and $N\sub M$ is a coisotropic submanifold then $\om|_N$ is a presymplectic structure on $N$ of corank equal to the codimension of $N$ in $M$. 

For a presymplectic manifold $(M,\om)$ we denote by $R^{M,\om}\sub M\x M$ its isotropy relation. By definition, this is the set of all pairs $(x(0),x(1))$, where $x\in C^\infty([0,1],M)$ is a path satisfying $\dot x(t)\in T_{x(t)}M^\om$, for every $t\in[0,1]$. For $x\in M$ we denote by $M^\om_x\sub M$ the isotropic leaf through $x$, i.e., the $R^{M,\om}$-equivalence class of $x$. 

We call $(M,\om)$ \emph{regular} if $R^{M,\om}$ is a closed subset and a submanifold of $M\x M$. Equivalently, there exists a smooth structure on the set of isotropic leaves $M_\om$ for which the canonical projection $\pi:M\to M_\om$ is a smooth submersion. In this case we define $\om_M$ to be the unique two-form on $M_\om$ such that $\pi^*\om_M=\om$. This is a symplectic form. 

For $i=0,1$ let $(M_i,\om_i)$ be a presymplectic manifold. We define the swap map 
\begin{equation}\label{eq:S M 1}S:M_1\x M_1\x M_2\x M_2\to M_1\x M_2\x M_1\x M_2\end{equation}
by $S\big(x_1,y_1,x_2,y_2\big):=\big(x_1,x_2,y_1,y_2\big)$.
\begin{lem}\label{le:R M 1 M 2} 
We have that $R^{M_1\x M_2,\om_1\oplus\om_2}=S\big(R^{M_1,\om_1}\x R^{M_2,\om_2}\big)$.
\end{lem}
\begin{proof}[of Lemma \ref{le:R M 1 M 2}]\setcounter{Claim}{0} This follows from a straight-forward argument.
\end{proof}
The next result is an immediate consequence of Lemma \ref{le:R M 1 M 2}. It is used in the proofs of Corollary \ref{cor:skinny}, Proposition \ref{prop:infinite}, and Lemma \ref{le:products} and Remark \ref{rmk:S N N'} below.
\begin{lem}\label{le:leaf prod} For every pair $(x_1,x_2)\in M_1\x M_2$, we have
\[(M_1\x M_2)^{\om_1\oplus\om_2}_{(x_1,x_2)}=(M_1)^{\om_1}_{x_1}\x(M_2)^{\om_2}_{x_2}.\]
\end{lem}
The next two results are used in the proof of Proposition \ref{prop:A Cross wt M}.
\begin{lem}\label{le:products}The identity map on $M_1\x M_2$ descends to a bijection 
\begin{equation}\label{eq:M 1 om 1 M 2 om 2}(M_1)_{\om_1}\x(M_2)_{\om_2}\to(M_1\x M_2)_{\om_1\oplus\om_2}.\end{equation}
If $(M_i,\om_i)$ is regular for $i=1,2$, then $\big(M_1\x M_2,\om_1\oplus\om_2\big)$ is regular and the map (\ref{eq:M 1 om 1 M 2 om 2}) is a symplectomorphism with respect to $(\om_1)_{M_1}\oplus(\om_2)_{M_2}$ and $(\om_1\oplus\om_2)_{M_1\x M_2}$. 
\end{lem}
\begin{proof}[of Lemma \ref{le:products}]\setcounter{Claim}{0} This follows from a straight-forward argument, using Lemmas \ref{le:R M 1 M 2} and \ref{le:leaf prod} and the definitions of the smooth and symplectic structures on the quotients.
\end{proof}
\begin{lem}\label{le:M 1 M 2} If the presymplectic manifold $\big(M_1\x M_2,\om_1\oplus\om_2\big)$ is regular then $(M_i,\om_i)$ is also regular, for $i=1,2$. 
\end{lem}
\begin{proof}[of Lemma \ref{le:M 1 M 2}]\setcounter{Claim}{0} It follows from Lemma \ref{le:R M 1 M 2} and Remark \ref{rmk:X i} below that $R^{N_i,\om_i}$ is a closed subset of $N_i\x N_i$, for $i=1,2$. Furthermore, Lemmas \ref{le:R M 1 M 2} and \ref{le:M i X i} imply that $R^{N_i,\om_i}$ is a submanifold of $N_i\x N_i$, for $i=1,2$. It follows that $(M_i,\om_i)$ is regular, for $i=1,2$. This proves Lemma \ref{le:M 1 M 2}.
\end{proof}

The next lemma is used in the proofs of Propositions \ref{prop:A Cross b} and \ref{prop:infinite}. Let $(M,\om)$ be a symplectic manifold. Recall the definitions (\ref{eq:S M om N},\ref{eq:S M om}).
\begin{lem}\label{le:S M om N S 2} Assume that $M$ is connected. Then the following statements hold:
\begin{enui}
\item\label{le:S M om N S 2:leq} For every non-empty coisotropic submanifold $N\sub M$ we have 
\begin{equation}\label{eq:S M om N S 2}S(M,\om,N)\cont S(M,\om).\end{equation}
\item\label{le:S M om N S 2:M} We have
\[S(M,\om,M)=S(M,\om).\]
\end{enui}
\end{lem}
For the proof of this lemma, we need the following result, which was also used in the proof of Theorem \ref{thm:A coiso cap}.
\begin{lem}\label{le:S M} We have
\begin{equation}\label{eq:S M om M sub}S(M,\om,M)\sub S(M,\om).\end{equation}
\end{lem}
\begin{proof}[of Lemma \ref{le:S M}] Let $u\in C^\infty(\D,M)$ be such that $u(S^1)$ is contained in some leaf of $M$. This leaf consists of a single point $x_0\in M$. We identify $S^2\iso\R^2\cup\{\infty\}$ and choose a map $f\in C^\infty(\D,S^2)$ that restricts to an orientation preserving diffeomorphism from $B^2$ to $\R^2$. We also choose a map $\rho\in C^\infty([0,1],[0,1])$ such that $\rho\const 1$ in a neighborhood of $0$ and $\rho(r)=1/r$ in a neighborhood of $1$. We define 
\[u':\D\to M,\quad u'(z):=u\big(\rho(|z|)z\big).\]
This map is constantly equal to $x_0$ in a neighborhood of $S^1$. Hence there exists a unique smooth map $v:S^2\to M$ satisfying $v\circ f=u'|_{B^2}$. We have
\[\int_{S^2}v^*\om=\int_\D{u'}^*\om=\int_\D u^*\om.\]
Here in the second equality we used the fact that $u'$ is smoothly homotopic to $u$ with fixed restriction to $S^1$. The inclusion (\ref{eq:S M om M sub}) follows. This proves Lemma \ref{le:S M}.
\end{proof}
\begin{proof}[of Lemma \ref{le:S M om N S 2}]\setcounter{Claim}{0} We prove {\bf statement (\ref{le:S M om N S 2:leq})}: Let $u\in C^\infty(S^2,M)$ be a map. We identify $S^2$ with $\R^2\cup\{\infty\}$. By Lemma \ref{le:homotop} below there exists a map $v\in C^\infty(S^2,M)$ that is smoothly homotopic to $u$ and satisfies $v(\infty)\in N$. 

We choose a smooth map $f:\D\to S^2$ that maps the interior $B_1\sub\D$ diffeomorphically and in an orientation preserving way onto $\R^2$. Then the map $v\circ f:\D\to M$ satisfies $v\circ f(S^1)\sub N$. Furthermore, 
\[\int_{S^2}u^*\om=\int_{S^2}v^*\om=\int_\D(v\circ f)^*\om.\]
The inclusion (\ref{eq:S M om N S 2}) follows. This proves (\ref{le:S M om N S 2:leq}).

{\bf Statement (\ref{le:S M om N S 2:M})} follows from statement (\ref{le:S M om N S 2:leq}) and Lemma \ref{le:S M}. This proves Lemma \ref{le:S M om N S 2}. 
\end{proof}
The next remark is used in the proofs of most main results of this paper. 
\begin{rem}\label{rmk:S N N'} Let $(M,\om)$ and $(M',\om')$ be symplectic manifolds, and $N\sub M$ and $N'\sub M'$ coisotropic submanifolds. Then 
\[S\big(M\x M',\om\oplus\om',N\x N'\big)=S(M,\om,N)+S(M',\om',N').\]
This follows from a straight-forward argument, using Lemma \ref{le:leaf prod}. $\Box$
\end{rem}
\begin{rem} Let $(M,\om)$ and $(M',\om')$ be symplectic manifolds of the same dimension, $N'\sub M'$ a coisotropic submanifold, and $\phi:M'\to M$ a symplectic embedding. The action spectrum $S(M',\om',N')$ is contained in $S(M,\om,\phi(N'))$. This follows from a straight-forward argument. $\Box$
\end{rem}
The next lemma gives a condition under which the opposite inclusion holds up to a correction term. It is used in the proofs of Theorem \ref{thm:A coiso cap}, Corollary \ref{cor:skinny}, and Propositions \ref{prop:infinite} and \ref{prop:ineqembA}. 
\begin{lem}\label{lemactionspec} If every continuous loop in a leaf of $N'$ is contractible in $M'$ then we have
\[S(M,\om,\phi(N'))\sub S(M',\om',N')+S(M,\om),\]
where the action spectrum $S(M,\om)$ is defined as in (\ref{eq:S M om}).
\end{lem}
\begin{proof}[of Lemma \ref{lemactionspec}]\setcounter{Claim}{0} Let $u\in C^{\infty}(\D,M)$ be a map such that $u(S^1)$ is contained in some isotropic leaf of $N:=\phi(N')$. It suffices to prove that
\begin{equation}\label{eq:int D S}\int_\D u^*\om\in S(M',\om',N')+S(M,\om).
\end{equation}
To see that this condition is satisfied, note that by our hypothesis the loop $x':=\phi^{-1}\circ (u|_{S^1}):S^1\to M'$ is contractible in $M'$. It follows that there exists a map $u'\in C^\infty(\D,M')$ such that $u'|_{S^1}=x'$. We denote by $\BAR{\D}$ the disk with the reversed orientation and by $\D\#\BAR{\D}$ the smooth oriented manifold obtained by concatenating the two disks along their boundary. We define $f:\D\#\BAR{\D}\to M$ to be the concatenation of $u$ and $\phi\circ u'$. It follows that 
\begin{equation}\label{eq:int D}\int_{\D\cup\BAR{\D}}f^*\om=\int_\D u^*\om-\int_\D(\phi\circ u')^*\om=\int_\D u^*\om-\int_\D{u'}^*\om'.\end{equation}
Since $\D\cup\BAR{\D}$ is diffeomorphic to $S^2$, we have $\int_{\D\cup\BAR{\D}}f^*\om\in S(M,\om)$. Combining this with (\ref{eq:int D}), the inclusion (\ref{eq:int D S}) follows. 
\end{proof}
The next result is used in the example on page \pageref{example A Cross} and the proofs of Corollaries \ref{cor:e w} and \ref{cor:skinny}, inequality (\ref{eq:A coiso d B}) in Theorem \ref{thm:A coiso cap}, and Theorem \ref{thm:emb d B}. For $k,n\in\N$ satisfying $k\leq n$ we denote by 
\[V(k,n):=\big\{\Theta\in\C^{k\x n}\,\big|\,\Theta\Theta^*=\one_k\big\}\]
the Stiefel manifold of unitary $k$-frames in $\C^n$.
\begin{prop}\label{minimalareaStiefel} We have
\[A(\C^{k\x n},\om_0,V(k,n))=\pi.\]
\end{prop}
\begin{proof}\setcounter{Claim}{0}
For a proof we refer to \cite[Proposition 1.3]{ZiLeafwise}.
\end{proof}
The next lemma is used in Remark \ref{rmk:phi L}. Recall the definitions (\ref{eq:Ham c M om},\ref{eq:Vert om c},\ref{eq:Vert phi}).
\begin{lem}\label{le:phi psi} Let $(M,\om)$ be a symplectic manifold, $K\sub M$ a compact subset, $\phi\in\Ham(M,\om)$, and $\eps>0$. Then there exists $\psi\in\Ham_c(M,\om)$ such that 
\begin{equation}\label{eq:psi K d c}\psi|_K=\phi|_K,\quad\Vert\psi\Vert_\om^c\leq\Vert\phi\Vert_\om+\eps.
\end{equation} 
(Here our convention is that $\infty+\eps:=\infty$.)
\end{lem}
For the proof of this lemma, we need the following. 
\begin{rem}\label{rem:flow}\upshape
Let $(M,\om)$ be a symplectic manifold, $H_0,H\in C^\infty(M,\R)$, $U\sub M$ an open subset, and $a>0$. Assume that $\phi_{H_0}^t$ (the Hamiltonian time $t$-flow of $H_0$) is well-defined on $U$ for $t\in[0,a]$, and that there exists a function $f\in C^\infty([0,1],\R)$ such that 
\[H\circ\phi_{H_0}^t=H_0\circ\phi_{H_0}^t+f(t)\]
on $U$. Then $\phi_H^t$ is well-defined on $U$ and $\phi_H^t=\phi_{H_0}^t$ on $U$, for $t\in[0,a]$. This follows from a straight-forward argument. $\Box$
\end{rem}
\begin{proof}[of Lemma \ref{le:phi psi}]\setcounter{Claim}{0} Without loss of generality we may assume that $M$ is connected and $\Vert\phi\Vert_\om<\infty$. We choose a function $H_0\in\HH(M,\om)$ (defined as on page \pageref{HH M om}) such that $\phi_{H_0}^1=\phi$ and $\Vert H_0\Vert\leq\Vert\phi\Vert_\om+\eps$. Since the set $K_0:=\bigcup_{t\in[0,1]}\phi_{H_0}^t(K)\sub M$ is compact, it has an open neighborhood $U_0\sub M$ with compact closure. We choose an open neighborhood $U_1\sub M$ of $\BAR{U_0}$ with compact closure. 

We choose a function $f\in C^\infty(M,[0,1])$ such that $f|_{M\wo U_1}\const0$ and $f|_{\BAR{U_0}}\const1$. We fix a point $x_0\in M$ and define
\[H:[0,1]\x M\to\R,\quad H(t,x):=f(x)\big(H_0(t,x)-H_0(t,x_0)\big).\]
Then the support of $H$ is contained in $\BAR{U_1}$ and hence compact. Furthermore, for $t\in[0,1]$ and $x\in U_0$ we have $H(t,\phi_{H_0}^t(x))=H_0(t,\phi_{H_0}^t(x))-H_0(t,x_0)$. Therefore, by Remark \ref{rem:flow} we have $\phi_H^1(x)=\phi_{H_0}^1(x)$, for every $x\in U_0$, and therefore the first condition in (\ref{eq:psi K d c}) holds. 

Finally, observe that 
\[\max_{x\in M}H(t,x)\leq\sup_{x\in M}H_0(t,x)-H_0(t,x_0).\]
Combining this with a similar inequality for $\min_{x\in M}H(t,x)$, it follows that $\Vert H\Vert\leq\Vert H_0\Vert$. Since $\Vert H_0\Vert\leq\Vert\phi\Vert_\om+\eps$, the second condition in (\ref{eq:psi K d c}) follows. This proves Lemma \ref{le:phi psi}.
\end{proof}
The next two remarks were used in the example on page \pageref{example A Cross} and the proofs of Theorem \ref{thm:A coiso cap}(\ref{thm:A coiso cap:A coiso d B Z}), Proposition \ref{prop:ineqembA}, and Corollary \ref{cor:exotic}. 
\begin{rem}\label{rem:e M M'} Let $(M,\om)$ and $(M',\om')$ be symplectic manifolds and $X\sub M$ a subset. Then we have
\[e\big(X\x M',M\x M'\big)\leq e(X,M).\]
This follows from a straight-forward argument. $\Box$
\end{rem}
\begin{rem}\label{rem:e Z 2n a} For every $n\in\N$ and $a>0$ we have
\[e(Z^{2n}(a),\R^{2n})\leq a.\]
This follows from a straight-forward argument. $\Box$
\end{rem}
\subsection{Topology and manifolds}\label{subsec:top mf}
The next result was used in the proof of Lemma \ref{le:S M om N S 2}.
\begin{lem}\label{le:homotop} Let $M$ and $M'$ be manifolds, $x_0\in M$, $x'_0\in M'$, and $u\in C^\infty(M',M)$. If $M$ is connected then $u$ is smoothly homotopic to a map $v\in C^\infty(M',M)$ satisfying $x_0=v(x'_0)$. 
\end{lem}
For the proof of this lemma we need the following two results. Let $n\in\N$, $M$ be a manifold, and $u\in C^\infty(\R^n,M)$. 
\begin{lem}\label{le:R n} There exists a map $h\in C^\infty([0,1]\x\R^n,M)$ satisfying 
\[h(1,x)=u(0),\quad\forall x\in B^n_1,\]
\begin{equation}\label{eq:h t x u}h(t,x)=u(x),\quad\forall (t,x)\in (\{0\}\x\R^n)\cup\big([0,1]\x(\R^n\wo B^n_2)\big).\end{equation}
\end{lem}
\begin{proof}[of Lemma \ref{le:R n}]\setcounter{Claim}{0} We choose a function $f\in C^\infty\big([0,1]\x[0,\infty),[0,\infty)\big)$ satisfying 
\begin{eqnarray*}&f(1,a)=0,\quad\forall a\in[0,1],&\\
&f(t,a)=1,\quad\forall(t,a)\in(\{0\}\x[0,\infty))\cup([0,1]\x[4,\infty)).&\end{eqnarray*}
We define $h(t,x):=u\left(f(t,|x|^2)x\right)$, for every $t\in[0,1]$ and $x\in\R^n$. This map has the required properties. This proves Lemma \ref{le:R n}.
\end{proof}
\begin{lem}\label{le:u ga} Let $\ga\in C^\infty([0,1],M)$ be a path. Assume that $u(x)=\ga(0)$, for every $x\in B^n_1$. Then there exists a map $h\in C^\infty([0,1]\x\R^n,M)$ satisfying $h(1,0)=\ga(1)$ and (\ref{eq:h t x u}).
\end{lem}
\begin{proof}[of Lemma \ref{le:u ga}]\setcounter{Claim}{0} We choose a function $f\in C^\infty\big([-1,1],[0,1]\big)$ satisfying $f(t)=0$ for $t\leq\frac12$, and $f(1)=1$. We define $h:[0,1]\x\R^n\to M$ by 
\[h(t,x):=\left\{\begin{array}{ll}\ga\circ f(t-|x|^2),&\textrm{if }|x|\leq1,\\
u(x),&\textrm{otherwise.}
\end{array}\right.\]
This map has the required properties. (Note that the hypothesis $u(x)=\ga(0)$, for every $x\in B_1^n$, ensures that $h$ is smooth along $[0,1]\x S^1$.) This proves Lemma \ref{le:u ga}.
\end{proof}
\begin{proof}[of Lemma \ref{le:homotop}]\setcounter{Claim}{0} Since $M$ is connected, there exists a path $\ga\in C^\infty([0,1],M)$ satisfying $\ga(0)=u(x'_0)$ and $\ga(1)=x_0$. The statement of the lemma follows from an argument using a chart around $x'_0$ and Lemmas \ref{le:R n} and \ref{le:u ga}. (We concatenate the homotopies provided by these Lemmas and smoothen the resulting homotopy.) This proves Lemma \ref{le:homotop}.
\end{proof}
The next two remarks and Lemma \ref{le:M i X i} below were used in the proof of Proposition \ref{prop:A Cross wt M}. These statements follow from elementary arguments.
\begin{rem}\label{rem:X Y} If $X_1$ and $X_2$ are non-empty topological spaces such that $X_1\x X_2$ is compact then $X_1$ and $X_2$ are compact. $\Box$
\end{rem}
\begin{rem}\label{rmk:X i} For $i=1,2$ let $X_i$ be a topological space and $\emptyset\neq A_i\sub X_i$ a subset. If $A_1\x A_2\sub X_1\x X_2$ is closed then $A_1$ and $A_2$ are closed. $\Box$
\end{rem}
\begin{lem}\label{le:M i X i} For $i=1,2$ let $M_i$ be a manifold and $X_i\sub M_i$ a non-empty subset. Assume that $X_1\x X_2\sub M_1\x M_2$ is a submanifold. Then $X_i$ is a submanifold of $M_i$, for $i=1,2$.
\end{lem}
The next lemma was used to ensure that the Hofer norm $\Vert\cdot\Vert_\om$ (see (\ref{eq:Vert phi})) is well-defined.
\begin{lem}\label{le:Borel} Let $X$ be a topological space and $f:[0,1]\x X\to \R$ a continuous function. Assume that there exists a sequence of compact subsets $K_\nu\sub X$, $\nu\in\N$ such that $\bigcup_\nu K_\nu=X$. Then the map 
\[[0,1]\ni t\mapsto \sup_{x\in X}f(t,x)\]
is Borel measurable. 
\end{lem}
\begin{proof} \setcounter{Claim}{0}This follows from an elementary argument. 
\end{proof}
\begin{acknowledgements} 
A considerable part of the work on this project was done during the second author's stay at the Max Planck Institute for Mathematics, Bonn. He would like to express his gratitude to the MPIM for the invitation and the generous fellowship. 
\end{acknowledgements}

\end{document}